\documentclass[preprint]{article}
\usepackage{amsfonts,amsmath,amssymb,amsthm,verbatim,placeins,color,comment,placeins,caption,subcaption,hyperref}
\usepackage{pgf}
\usepackage{tikz}
\usetikzlibrary{arrows,automata,positioning}
\usepackage{csquotes}
\usepackage{graphicx}
\usepackage{epstopdf}
\usepackage{url}
\usepackage{float}
\usepackage{array}
\usepackage[english]{babel}
\usepackage[T1]{fontenc}
\usepackage[utf8]{inputenc}
\usepackage[title]{appendix}

\newtheorem{definition}{Definition}[section]
\newtheorem{remark}{Remark}[section]
\newtheorem{exa}{Example}[section]
\newtheorem{theorem}{Theorem}[section]

\tikzset{
	style1/.style = {fill=lightgray,draw=black,text=black},
	style2/.style = {fill=white,draw=none,text=white}
}

\begin{document}
	
	\title{\textbf{Epidemic random model under vaccination regime and non-permanent immunity driven by time-changed L\'{e}vy noise}}
	\maketitle
	
	\centerline{
		\author{\textbf{Giulia Di Nunno}$^{\textrm{a, b}}$,
			\textbf{Jasmina \DJ or\dj evi\'{c}}$^{\textrm{c}}$,
			\textbf{Nenad \v{S}uvak}$^{\textrm{d}}$
	}}
	
	\bigskip \bigskip

	$^{\textrm{a}}$\emph{Department of Mathematics, University of Oslo, Blindern 0316 Oslo, Norway}\footnote{e-mail: giulian@math.uio.no}

	$^{\textrm{b}}$\emph{Department of Business and Management Science, NHH Norwegian School of Economics, \\ \phantom{aaaa}Bergen, Norway}

	$^{\textrm{c}}$\emph{Faculty of Sciences and Mathematics, University of Ni\v s,  Vi\v{s}egradska 33, 18000 Ni\v{s}, Serbia}\footnote{e-mail: jasmina.djordjevic@pmf.edu.rs}

	$^{\textrm{d}}$\emph{J. J. Strossmayer University of Osijek, School of Applied Mathematics and Informatics, \\ \phantom{aaaa}Trg Ljudevita Gaja 6, 31000 Osijek, Croatia}\footnote{e-mail: nsuvak@mathos.hr}
	
	\bigskip \bigskip
	
	\noindent\textbf{Abstract.}
	A random version of the Susceptible-Infected-Recovered-Vaccinated (SIRV) epidemic model with non-permanent immunity and a population of non-constant size is studied. The compartments satisfy a system of differential equations in which the contact rate is a random coefficient, and the contact rate itself is driven by a stochastic differential equation with time-changed L\'{e}vy noise. This noise is constructed from a conditional Brownian motion and a doubly stochastic Poisson random field, both governed by a stochastic clock -- when the clock is independent of the Brownian and Poisson components, the two parts of the noise are a time-changed Brownian motion and a time-changed Poisson random measure. The equation for the contact rate decouples from the compartment equations and enters it as a random parameter, which allows the compartments to be treated pathwise. The existence, uniqueness and positivity of the global solution is proven and the basic reproduction number of the model is shown to be random. Furthermore, sufficient conditions for extinction and for persistence in mean of the disease are derived. The conditions are expressed through the parameters of the model and two long-run averages of the noise intensity, which are governed by the average speed of the stochastic clock.

	\bigskip


	\bigskip

	
	\bigskip

	\section{Introduction} \label{intro}

	A random version of the Susceptible-Infected-Recovered-Vaccinated (SIRV) epidemic model with non-permanent immunity is studied in this paper. The starting point is the deterministic SIRV model introduced in  \cite{im}, which is modified to allow a population of non-constant size, bounded by an environmental capacity $K$, with natural and disease-induced death rates. Randomness enters the model through the contact rate modeled by a stochastic process $\widetilde{\beta} = ( \widetilde{\beta}(t), \, t \geq 0)$ which is driven by a time-changed L\'{e}vy noise. The resulting model is a partially coupled system of differential equations, with four random differential equations (RDEs) for the compartments, and one stochastic differential equation (SDE) for the contact rate, which enters the first four dynamics as a random coefficient.

	For this model we prove the existence, uniqueness and positivity of the global solution, derive the basic reproduction number by the next-generation method and show that, being a function of $\widetilde{\beta}$, it is itself a stochastic process. We also derive sufficient conditions for the extinction and for the persistence in mean of the disease. Both thresholds are expressed through the parameters of the model and the average speed of the stochastic clock. To the best of our knowledge, models of this type with time-changed noises have not yet been studied in the epidemic literature.

	The need to extend deterministic epidemic models by stochastic components has long been recognised (see e.g. \cite{Allen2008}-\cite{1n}, \cite{Gray2021,Mao2002, Pardoux2015}), and the COVID-19 pandemic has shown the importance of models able to capture uncertainty, variability and abrupt changes in transmission. 
    To capture uncertainty, the simplest and most widely used approach is the additive perturbation of a constant contact rate, $\beta \, dt \mapsto \beta \, dt + dB(t)$, i.e. the case of a constant drift and a constant diffusion coefficient, possibly enriched by a pure-jump or L\'{e}vy component (see e.g. \cite{moi1}-\cite{ moi4}, \cite{9n, 8n}). Another approach is to model the contact rate by an It\^{o} process,
	$$d\beta(t) = \mu(t) \, dt + \sigma(t) \, dB(t), \quad t \geq 0,$$
	with drift and diffusion coefficients possibly depending on $\beta(t)$ (see e.g. \cite{Wang2018}, where an Ornstein--Uhlenbeck process is incorporated in the SIS model). The contact rate can also be modeled by an SDE with L\'{e}vy jumps in addition to the Brownian component (see e.g. \cite{Zhang2017}), where the jumps represent abrupt environmental changes. These systems are analysed with a fixed set of tools: the existence and uniqueness of a positive local solution follow from the classical theory (see e.g. \cite{mao, OS-book}), the extension to a global solution relies on showing that the explosion time is a.s. infinite, and the explicit conditions for extinction and persistence follow from the multidimensional It\^{o} formula applied to a suitably chosen Lyapunov function (see e.g. \cite{Korobeinikov2002}), together with elementary stochastic estimates. The same tools are used here, with the difference that the compartments themselves carry no noise. 
	
	Memory can be introduced into the dynamics through delays and Volterra-type equations (see e.g. \cite{7n, 8nn} and the references therein). Another way of introducing memory into a stochastic model is a time change, since the operational time flows irregularly and the increments of the noise are no longer independent of its past. Time-changed processes replace the calendar time by a stochastic clock and thereby account for irregular progression of events, bursts and slow-downs of the epidemic and heavy-tailed waiting times. Additionally, a time change makes it possible to include the clustering of transmission which is to be expected as a consequence of social behaviour, in particular of high concentrations of individuals in mostly indoor settings, such as gatherings, public transport, restaurants, workplaces or care homes (see e.g. \cite{Leclerc2020}, where $201$ clusters of SARS-CoV-2 transmission are classified by the setting). 
    
    In this paper the time change acts on the noise which drives the contact rate. This noise is the martingale random field $\mu = B + \widetilde{H}$ introduced and studied in \cite{GS1, GS2}, i.e. the sum of a conditionally Gaussian random measure $B$, accounting for continuous fluctuations, and a centred doubly stochastic Poisson random measure $\widetilde{H}$, accounting for jumps. Both components are governed by a random measure $\Lambda$, the stochastic clock, whose intensity is given by two non-negative processes $\lambda^{B}$ and $\lambda^{H}$. The clock modulates the variance of the continuous fluctuations and the frequency of the jumps. It may capture seasonality, interventions such as lockdowns or vaccination campaigns and the time clustering of infections. Conditionally on the clock, $B$ and $\widetilde{H}$ are a time-changed Brownian motion and a time-changed centered Poisson random measure, and therefore the stochastic calculus of \cite{GS2} is available. In this framework, the intensities of the clock must not depend on the noise they govern, as is the case for deterministic periodic intensities, Markov-modulated intensities, positive diffusions such as the Cox--Ingersoll--Ross process, or L\'{e}vy-driven Ornstein--Uhlenbeck intensities. Self-exciting intensities of Hawkes type, in which $\lambda^{H}$ is a functional of the past of the jump measure itself, fall outside this framework, since conditioning on the whole path of the intensity then already determines the jump times (for a stochastic epidemic model with Hawkes-type self-excitation we refer to \cite{Tymoshenko-2026}).

	In stochastic compartmental models the noise usually perturbs the compartment equations directly, so that the compartments themselves are semimartingales. For example, in the stochastic SIRV model studied in \cite{im} the contact rate is perturbed by a Gaussian white noise, $\beta \mapsto \beta + \varepsilon \eta(t)$, so that the incidence terms (and therefore the equations for the susceptible, infected and vaccinated classes) contain Brownian components, while in \cite{Tymoshenko-2025, Tymoshenko-2026} every compartment equation carries its own noise. In this paper, the
compartments satisfy RDEs and the only random object is the contact rate $\widetilde{\beta}$ through which behaviour,
policy and environment act on an epidemic. The contact rate is also the parameter which is hardest to specify -- modeling it as a stochastic process with its own dynamics, driven by a noise whose clock is itself stochastic, keeps the mechanism of transmission deterministic given the contact rate. Technically, the compartment RDEs can then be treated pathwise, while the stochastic analysis is confined to the scalar equation for the contact rate. Specifically, we model the contact rate as
	\begin{equation*}
		\widetilde{\beta}(t) = \widetilde{\beta}_{max} f\left( \beta(t) \right), \qquad d\beta(t) = \beta(t) \, dt + \int\limits_{\mathbb{R}} \sigma_{t}(z, \beta(t)) \, \mu(dt, dz),
	\end{equation*}
	where $f$ is a smooth function with values in $[f_{min}, 1]$, $0 < f_{min} < 1$, so that $\widetilde{\beta}$ takes values between $\widetilde{\beta}_{max} f_{min}$ and $\widetilde{\beta}_{max}$. The positive constant $\widetilde{\beta}_{max}$ is interpreted as the ceiling of the contact rate, i.e. the largest average number of adequate contacts per individual per unit of time, and $\sigma_{t}(z, b)$ determines the impact of a shock of size $z$ on the underlying process $\beta$. The identity drift in the equation for $\beta$ is the simplest representative of an affine drift. A general affine drift, including the mean-reverting case, would not change the arguments but would burden the notation (see Remark \ref{rem-affine-drift}).

	The paper is organised as follows. Section~\ref{subsec-preliminaries} recalls the construction and the main properties of the time-changed L\'{e}vy noise. Section~\ref{sec:model} introduces the model: the deterministic SIRV system, the stochastic contact rate, the random SIRV model and its basic reproduction number. Section~\ref{sec:EUPsolution} establishes the existence, uniqueness and positivity of the global solution and identifies the positively invariant set of the model. Sections~\ref{sec:extinction} and~\ref{sec:persistence} give the conditions for the extinction and for the persistence in mean of the disease. The model and these results are illustrated by simulations throughout the paper, and Section~\ref{sec:simulation} gives the parameter values, based on the data and published estimates for SARS-CoV-2 in Denmark, together with the simulation setting. Section~\ref{sec:conclusion} contains some concluding remarks. The construction of the stochastic integral with respect to the driving noise is recalled in Appendix \ref{app:noise}.

	\section{Preliminaries} \label{subsec-preliminaries}

	In this section we recall the properties of the random measure $\mu$, acting as a driving noise in the contact rate model. In particular, we address its construction from a stochastic clock and its martingale property. These results are due to \cite{GS2}, where further properties of this random measure can be found.

	Let $(\Omega, \mathcal{F}, P)$ be a complete probability space. Further, let $\mathcal{X} := [0, \infty \rangle \times \mathbb{R} = \left([0, \infty \rangle \times \{0\}\right) \cup \left([0, \infty \rangle \times \mathbb{R}_{0}\right)$, where $\mathbb{R}_{0} := \mathbb{R} \setminus \{0\}$, let $\mathcal{B}_{\mathcal{X}}$ be the Borel $\sigma$-algebra on $\mathcal{X}$ and let $\Delta$ denote a Borel set from $\mathcal{B}_{\mathcal{X}}$. Each component $\lambda^{l}$, $l = B, H$
	of the two-component stochastic process  $\lambda := \left(\lambda^{B}, \lambda^{H}\right)$  has the following properties:	
	\begin{itemize}
		\item[\it{(l1)}] $\lambda^{l}$ is non-negative $P$-a.s., i.e. $\lambda_{t}^{l} \geq 0$ $P$-a.s. for all $t \geq 0$,
		\item[\it{(l2)}] $\lambda^{l}$ is stochastically continuous, i.e. $\lim\limits_{h \rightarrow 0} P \left( \left| \lambda_{t+h}^{l} - \lambda_{t}^{l} \right| \geq \varepsilon \right) = 0$ for all $\varepsilon > 0$ and all $t \geq 0$,
		\item[\it{(l3)}] $E \left[ \int_{0}^{t} \lambda_{s}^{l} \, ds \right] < \infty$, and therefore $\int_{0}^{t} \lambda_{s}^{l} \, ds < \infty$ $P$-a.s., for every $t \geq 0$.
	\end{itemize}
	Conditions (l1)-(l3) are needed to define the stochastic clock that runs with intensity (rate) determined by the process $\lambda = \left(\lambda^{B}, \lambda^{H}\right)$. Specifically, (l3) is required on bounded time-intervals only, so that the stochastic clock is allowed to keep running at a positive average speed indefinitely, which is needed for the long-run thresholds derived in Sections \ref{sec:extinction} (extinction of the disease) and \ref{sec:persistence} (persistence of the disease).

	This particular stochastic clock, determining the operational time of the random measure $\mu$, is a random measure $\Lambda$ on $\mathcal{B}_{\mathcal{X}}$ defined in the following way:
	\begin{equation*} \label{Lambda-defininition}
		\Lambda(\Delta) := \int_{0}^{\infty} \mathbf{1}_{\{(t, 0) \in \Delta\}}(t) \lambda_{t}^{B} \, dt + \int_{0}^{\infty} \int_{\mathbb{R}_{0}}  \mathbf{1}_{\Delta}(t, z) \, \nu(dz) \lambda_{t}^{H} \, dt,
	\end{equation*}
	where $\nu$ is a $\sigma$-finite measure on the Borel subsets of $\mathbb{R}_{0}$, satisfying
	\begin{equation} \label{measure-nu}
		\int_{\mathbb{R}_{0}} z^{2} \, \nu(dz) < \infty.
	\end{equation}
	Next, we write $\Lambda^{B}$ and $\Lambda^{H}$ for the restrictions of $\Lambda$ to $[0, \infty \rangle \times\{0\}$ and to $[0, \infty \rangle \times \mathbb{R}_{0}$, respectively, so that
	\begin{equation*}  \label{Lambda-representation}
		\Lambda \left(\Delta \right) = \Lambda^{B} \left( \Delta \cap [0, \infty \rangle \times\{0\} \right) + \Lambda^{H} \left( \Delta \cap [0, \infty \rangle \times \mathbb{R}_{0} \right), \quad \Delta \in \mathcal{B}_{\mathcal{X}},
	\end{equation*}
	with $\Lambda^{B} \left( ds, \{0\} \right) = \lambda^B_{s} \, ds$ and $\Lambda^{H} \left( ds, dz \right) = \lambda^H_{s} \, ds \, \nu(dz)$ for $z \in \mathbb{R}_{0}$, and we let $\mathcal{F}^{\Lambda}_{\infty}$ be the $\sigma$-algebra generated by the values of $\Lambda$. 
	
	In the epidemiological context of Section \ref{sec:model}, $\Lambda$ is the stochastic clock of the epidemic, while $\lambda^{B}$ and $\lambda^{H}$ are the instantaneous speeds at which the operational time of the continuous fluctuations and of the shocks of the contact rate \eqref{SDE-beta} accumulate with respect to the calendar time. For better understanding of the nature of the time-change in the definition of the random measure $\mu$, the following two random measures are defined.

	\begin{definition} \label{B-H-definition}
		$B$ is a signed random measure on the Borel sets of $[0, \infty \rangle \times \{0\}$ satisfying:
		\begin{itemize}
			\item[(A1)] $P \left( B(\Delta) \leq x \mid \mathcal{F}^{\Lambda}_{\infty} \right) = P \left( B(\Delta) \leq x \mid \Lambda^{B}(\Delta) \right) = \Phi \left( \frac{x}{\sqrt{\Lambda^{B}(\Delta)}}\right), x \in \mathbb{R}, \Delta \subseteq [0, \infty \rangle \times \{0\}$,
			\item[(A2)] $B \left( \Delta_{1} \right)$ and $B\left( \Delta_{2} \right)$,  whenever $\Delta_{1}$ and $\Delta_{2}$ are disjoint, are conditionally independent given $\mathcal{F}^{\Lambda}_{\infty}$.
		\end{itemize}
		Here $\Phi$ is the cumulative distribution function of a standard normal random variable. $H$ is a random measure on the Borel sets of $[0, \infty \rangle \times \mathbb{R}_{0}$ satisfying:
		\begin{itemize}
			\item[(A3)] $P \left( H(\Delta) = k \mid \mathcal{F}^{\Lambda}_{\infty} \right) = P\left( H(\Delta) = k \mid \Lambda^{H}(\Delta) \right) = \frac{\Lambda^{H}(\Delta)^{k}}{k !} e^{-\Lambda^{H}(\Delta)}, k \in \mathbb{N}_{0}, \Delta \subseteq [0, \infty \rangle \times \mathbb{R}_{0}$,
			\item[(A4)] $H \left( \Delta_{1} \right)$ and $H \left( \Delta_{2} \right)$, whenever $\Delta_{1}$ and $\Delta_{2}$ are disjoint, are conditionally independent given $\mathcal{F}^{\Lambda}_{\infty}$.
		\end{itemize}
		Furthermore, it is assumed that
		\begin{itemize}
			\item[(A5)] $B$ and $H$ are conditionally independent given $\mathcal{F}^{\Lambda}_{\infty}$.
		\end{itemize}
	\end{definition}
	
	In (A1) and (A3) the sets $\Delta$ are such that $\Lambda(\Delta) < \infty$ $P$-a.s., which is what makes the conditional distributions well defined. By (l3) this covers all sets bounded in time for $B$, and all sets of the form $[0, t] \times A$ with $\nu(A) < \infty$ for $H$. 
	Further, let $\widetilde{H} := H - \Lambda^{H}$ be the signed random measure given by $\widetilde{H}(\Delta) = H(\Delta) - \Lambda^{H}(\Delta)$, $\Delta \subset [0, \infty \rangle \times \mathbb{R}_{0}$, $\Lambda^{H}(\Delta) < \infty$. Having all these random measures at hand, we can define the random measure $\mu$.

	\begin{definition} \label{mu-definition}
		The signed random measure $\mu$ on the Borel subsets of $\mathcal{X}$ is defined by
		\begin{equation*}
			\mu(\Delta) := B \left( \Delta \cap [0, \infty \rangle \times \{0\} \right) + \widetilde{H} \left( \Delta \cap [0, \infty \rangle \times \mathbb{R}_{0} \right), \quad \Delta \in \mathcal{B}_{\mathcal{X}}.
		\end{equation*}
	\end{definition}

	From (A1) and (A3) from Definition \ref{B-H-definition}, conditionally on $\mathcal{F}^{\Lambda}_{\infty}$, the first moments of $B$ and $\widetilde{H}$ vanish and their second moments are $E \left[ B(\Delta)^{2} \mid \mathcal{F}^{\Lambda}_{\infty} \right] = \Lambda^{B}(\Delta)$ and $E \left[ \widetilde{H}(\Delta)^{2} \mid \mathcal{F}^{\Lambda}_{\infty} \right] = \Lambda^{H}(\Delta)$. Therefore, by the conditional independence (A2), (A4) and (A5) from Definition \ref{B-H-definition}, it follows that
	\begin{equation} \label{mu-conditional-moments}
		E \left[ \mu(\Delta) \mid \mathcal{F}^{\Lambda}_{\infty} \right] = 0, \qquad
		E \left[ \mu(\Delta)^{2} \mid \mathcal{F}^{\Lambda}_{\infty} \right] = \Lambda(\Delta), \qquad
		E \left[ \mu \left( \Delta_{1} \right) \mu \left( \Delta_{2} \right) \mid \mathcal{F}^{\Lambda}_{\infty} \right] = 0,
	\end{equation}
	whenever  are disjoint. Therefore, the values of the random measure $\mu$ on disjoint sets are conditionally orthogonal, and $\mu$ is additive.

	The term time-change related to the measure $\Lambda$ is justified by the following characterisation, due to \cite{Serfozo} (see also \cite{BR}). Denote $B_{t} := B([0, t] \times \{0\})$, $\Lambda_{t}^{B} := \int_{0}^{t} \lambda_{s}^{B} \, ds$, $\eta_{t} := \int_{0}^{t} \int_{\mathbb{R}_{0}} z \, \widetilde{H}(ds, dz)$ and $\widehat{\Lambda}_{t}^{H} := \int_{0}^{t} \lambda_{s}^{H} \, ds$.

	\begin{theorem} \label{theorem-time-change}
		Let $W = \left( W_{t}, t \geq 0 \right)$ be a Brownian motion and $N = \left( N_{t}, t \geq 0 \right)$ be a centred pure jump L\'{e}vy process with L\'{e}vy measure $\nu$. Assume that both $W$ and $N$ are independent of $\Lambda$. Then $B$ satisfies (A1) and (A2) if and only if, for any $t \geq 0$, $B_{t} \stackrel{d}{=} W_{\Lambda_{t}^{B}}$, and $\eta$ satisfies (A3) and (A4) if and only if, for any $t \geq 0$, $\eta_{t} \stackrel{d}{=} N_{\widehat{\Lambda}_{t}^{H}}$.
	\end{theorem}

	Conditionally on $\mathcal{F}^{\Lambda}_{\infty}$, $B$ is thus a Brownian motion and $\eta$ a centred pure jump L\'{e}vy process, running in the operational times $\Lambda^{B}_{t}$ and $\widehat{\Lambda}^{H}_{t}$, respectively. For the present paper the following two statements of Theorem \ref{theorem-time-change} are of interest:
	\begin{itemize}
		\item if $B$ satisfies (A1) and (A2) from the Definition \ref{B-H-definition}, there exists a Brownian motion $W$, independent of $\Lambda$, such that $B_{t} \stackrel{d}{=} W_{\Lambda_{t}^{B}}$ for every $t \geq 0$,
		\item analogously, if $\eta$ satisfies (A3) and (A4) from the Definition \ref{B-H-definition}, there exists a centred pure jump L\'{e}vy process $N$ with L\'{e}vy measure $\nu$, independent of $\Lambda$, such that $\eta_{t} \stackrel{d}{=} N_{\widehat{\Lambda}_{t}^{H}}$ for every $t \geq 0$.
	\end{itemize}
	
	\begin{remark} \label{rem-Lambda-generalisation}
	\begin{itemize}
		\item[]
		\item[(a)] The assumption \eqref{measure-nu} on the measure $\nu$ is more strict than the integrability of $(1 \wedge z^{2})$ required for a general L\'{e}vy measure, but it ensures that the jump coefficient of the contact rate model \eqref{SDE-beta} is square-integrable with respect to the measure $\Lambda$. No finiteness of $\nu(\mathbb{R}_{0})$ is assumed, and therefore infinitely many small jumps on bounded time-intervals are allowed. 

		\item[(b)] By Theorem \ref{theorem-time-change}, the processes $(B_{t}, \, t \geq 0)$ and $(\eta_{t}, \, t \geq 0)$ are L\'{e}vy processes on a stochastic time-scale, so it is worth mentioning in which sense they may fail to be L\'{e}vy processes in the calendar time. A L\'{e}vy process is characterised by stationary and independent increments, and the stochastic clock $\Lambda$ influences both of these properties.
		
		Stationarity of increments is lost as soon as stochastic clock does not run at a constant rate with respect to the calendar time. In particular, conditionally on $\mathcal{F}^{\Lambda}_{\infty}$, the increment of the random measure $B$ over $\langle s, t]$ is centred Gaussian with variance $\Lambda^{B} \left( \langle s, t] \times \{ 0 \} \right) = \int_{s}^{t} \lambda^{B}_{u} \, du$, which depends on $(t - s)$ alone if and only if $\lambda^{B}$ is time-independent.
		
		Independence of increments is lost as soon as the clock is stochastic. In that case, the increment carries information from $\mathcal{F}^{\Lambda}_{\infty}$, and therefore about the distribution of the subsequent increments. 
		
		What survives the introduction of the stochastic clock is the conditional independence of increments given $\mathcal{F}^{\Lambda}_{\infty}$ stated in (A2) and (A4) in Definition \ref{B-H-definition}, and this is the exact property that characterises the class of noises considered in this paper (see e.g. \cite{BR}).
		
		\item[(c)] The simplest example of the stochastic clock $\Lambda$ is obtained for constant intensities, $\lambda^{B} = c_{B}$ and $\lambda^{H} = c_{H}$. In this case $\Lambda(dt, dz) = \left( c_{B} \mathbf{1}_{\{ z = 0 \}} + c_{H} \mathbf{1}_{\{ z \neq 0 \}} \nu(dz) \right) dt$ is deterministic. For deterministic, but non-constant intensities $\lambda^{B}$ and $\lambda^{H}$, the increments remain independent but not stationary, and the corresponding processes are just additive.
		
		Furthermore, by (A3) from Definition \ref{B-H-definition}, $H$ is a doubly stochastic Poisson random measure, i.e. a Cox random measure governed by $\Lambda^{H}$. The admissible intensities are those which do not depend on the noise they govern, such as those given in Example \ref{ex-clocks} below. Self-exciting intensities of Hawkes type, in which $\lambda^{H}$ is a functional of the past of $H$ itself, are not admissible, since conditioning on the whole path of the intensity would already determine the jump times and (A3) would fail (see also Section \ref{sec:conclusion}).
	\end{itemize}
	\end{remark}

	Further, let $\mathbb{F}^{B} = \left( \mathcal{F}_{t}^{B}, t \geq 0 \right)$ and $\mathbb{F}^{H} = \left( \mathcal{F}_{t}^{H}, t \geq 0 \right)$ be the filtrations generated by $B$ and by $H$, respectively, let $\mathbb{F}^{\Lambda} = \left( \mathcal{F}_{t}^{\Lambda}, t \geq 0 \right)$ be the filtration generated by $\Lambda$, and let $\mathcal{F}^{\mu}_{t} := \mathcal{F}_{t}^{B} \vee \mathcal{F}_{t}^{H} \vee \mathcal{F}_{t}^{\Lambda}$ be the $\sigma$-algebra generated by the random measure $\mu$ and by the stochastic clock $\Lambda$ up to time $t$. The filtration used in this paper is the right-continuous version $\mathbb{F} = \left( \mathcal{F}_{t}, t \geq 0 \right)$,
	\begin{equation}  \label{F-filtration}
		\mathcal{F}_{t} = \bigcap_{s > t} \mathcal{F}_{s}^{\mu},
	\end{equation}
	to which $\mu$ is adapted. For $\Delta \subset [t, \infty \rangle \times \mathbb{R}$ with $\Lambda(\Delta) < \infty$ $P$-a.s., according to \it{(A1)--(A4)} from the Definition \ref{B-H-definition} it follows:
	\begin{equation*}
		E \left[ \mu(\Delta) \mid \mathcal{F}_{t} \right] = E \left[ E \left[ \mu(\Delta) \mid \mathcal{F}_{t} \vee \mathcal{F}^{\Lambda}_{\infty} \right] \mid \mathcal{F}_{t} \right] = 0,
	\end{equation*}
	i.e. $\mu$ has the martingale property with respect to $\mathbb{F}$. Together with its additivity and the conditional orthogonality \eqref{mu-conditional-moments}, it means that $\mu$ is a martingale random field with respect to $\mathbb{F}$ (see \cite{GS2}, Remark $2.6$).

	\begin{exa}[Stochastic clocks] \label{ex-clocks}
		The following intensities satisfy requirements (l1)-(l3), and also (l4) additionally introduced in Section \ref{sec:extinction} ($\ell = \lim_{t \to \infty} \frac{1}{t} \int_{0}^{t} \lambda_{s} \, ds < \infty$), and either component of the process $\lambda = (\lambda^{B}, \lambda^{H})$ may be chosen from these examples:
		\begin{itemize}
			\item[(i)] Constant intensity $\lambda_{t} = c > 0$ (yields the classical L\'{e}vy case),
			\item[(ii)] Deterministic periodic intensity $\lambda_{t} = c \left( 1 + a \sin (2 \pi t / p) \right)$, $c > 0$, $|a| < 1$, $p > 0$,
			\item[(iii)] Markov-modulated intensity $\lambda_{t} = \lambda_{J(t)}$, where $J$ is an irreducible continuous-time Markov chain on $\{1, \ldots, m\}$ with stationary distribution $\pi$, independent of $B$ and $H$, and $\lambda_{1}, \ldots, \lambda_{m} > 0$,
			\item[(iv)] Cox-Ingersoll-Ross intensity $d \lambda_{t} = a \left( b - \lambda_{t} \right) dt + \varsigma \sqrt{\lambda_{t}} \, d W_{t}$, $a, b, \varsigma > 0$, $2 a b \geq \varsigma^{2}$, where Brownian motion $W$ is independent of random measures $B$ and $H$,
			\item[(v)] L\'{e}vy-driven Ornstein-Uhlenbeck intensity $d \lambda_{t} = - a \lambda_{t} \, dt + d L_{a t}$, $a > 0$, where a subordinator $L$ is independent of $B$ and $H$ and $E[L_{1}] < \infty$.
		\end{itemize}
	\end{exa}

\begin{figure}[htbp]
	\centering
	\includegraphics[width=\textwidth]{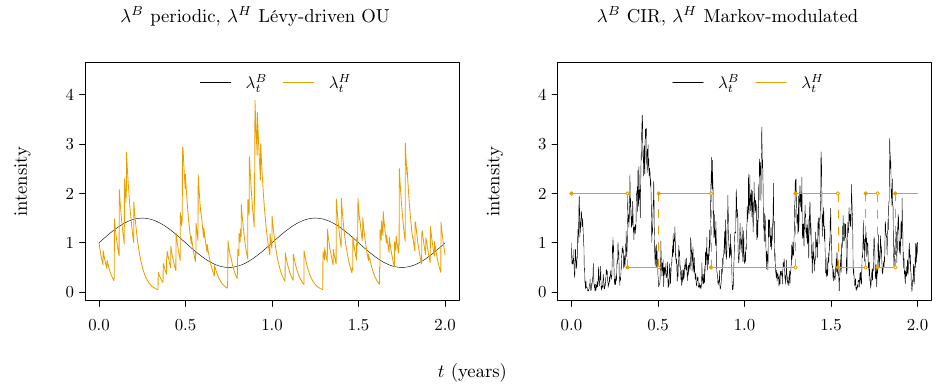}
	\caption{Sample paths of the processes $\lambda^{B}$ and $\lambda^{H}$, for the periodic $\lambda^{B}$ and the L\'{e}vy-driven Ornstein--Uhlenbeck $\lambda^{H}$ (left) and for the Cox--Ingersoll--Ross $\lambda^{B}$ and the Markov-modulated $\lambda^{H}$ (right), with parameter values given in Section~\ref{sec:simulation}.} \label{fig:clocks-intensity}
\end{figure}

The stochastic integral with respect to $\mu$, which gives a meaning to the contact rate model \eqref{SDE-beta}, is constructed in \cite{GS2}, Section $2$. Its construction, the isometry it satisfies and the quadratic variations of the martingale it produces are recalled in Appendix \ref{app:noise}, together with the facts about quadratic variations of general semimartingales used in the proofs of results in this paper.

\begin{figure}[H]
	\centering
	\includegraphics[width=\textwidth]{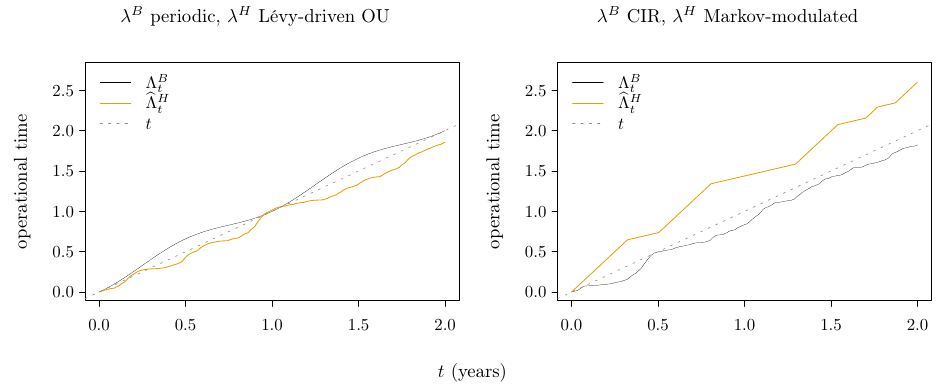}
	\caption{Sample paths of the operational times $\Lambda^{B}_{t}$ and $\widehat{\Lambda}^{H}_{t}$ corresponding to the processes $\lambda^{B}$ and $\lambda^{H}$ from Figure~\ref{fig:clocks-intensity}, with parameter values given in Section~\ref{sec:simulation}.} \label{fig:clocks-optime}
\end{figure}

\begin{figure}[H]
	\centering
	\includegraphics[width=\textwidth]{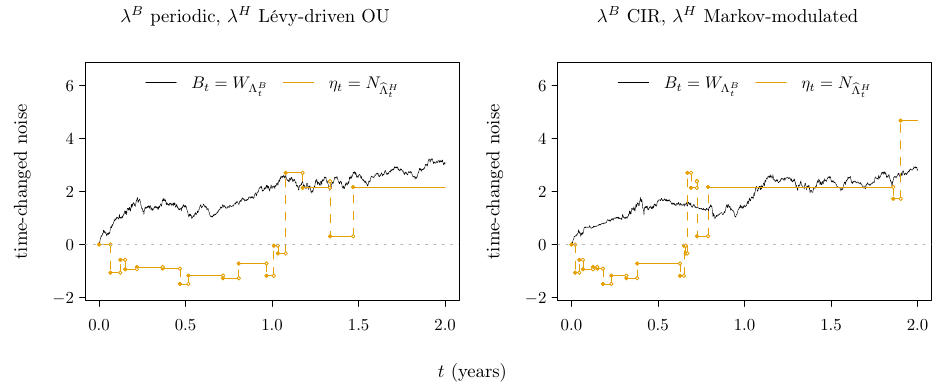}
	\caption{Sample paths of the time-changed processes $B_{t} = W_{\Lambda^{B}_{t}}$ and $\eta_{t} = N_{\widehat{\Lambda}^{H}_{t}}$ corresponding to the operational times from Figure~\ref{fig:clocks-optime}, with parameter values given in Section~\ref{sec:simulation}. } \label{fig:clocks-noise}
\end{figure}

\section{SIRV model driven by time-changed L\'{e}vy noise} \label{sec:model}
	
In general, mathematical model of an epidemic describes the dynamics of the spread of infectious disease within a population. In compartmental models the population is divided into finitely many disjoint compartments or classes, based on the disease status of individuals, and the model describes the transitions of individuals between the compartments over time. Such models are useful for better understanding of the epidemic dynamics, for predicting the spread of the disease and for evaluating the impact of interventions such as vaccination, social distancing and quarantine.

\subsection{Deterministic SIRV model} \label{SIRV-deterministic}
	
The compartmental epidemic model presented in this section is inspired by the deterministic SIRV model studied in \cite{im}, in which the population is divided into four disjoint compartments of individuals: $S$ (susceptible), $I$ (infected), $R$ (recovered) and $V$ (vaccinated). In this particular model, the total size of the population at time $t$, $N(t) = S(t) + I(t) + R(t) + V(t)$, is constant due to the equality of the birth and the death rates. We retain the four compartments but let the population size vary, and we distinguish the natural from the disease-induced death rate. Our model depends on the following parameters:
	\begin{itemize}
		\item $\beta \geq 0$ is the contact rate, i.e. the average number of adequate contacts of one individual per unit of time (an adequate contact is the contact which is sufficient for the transmission of the infection),
		\item $\lambda > 0$ is the birth rate,
		\item $\kappa > 0$ is the natural death rate,
		\item $\kappa_1 \geq \kappa$ is the disease-induced death rate,
		\item $\rho > 0$ is the vaccination rate of susceptible individuals,
		\item $\delta \in [0, 1]$ is the vaccine efficiency parameter ($\delta = 1$ -- useless vaccine, $\delta=0$ -- total protection),
		\item $\theta > 0$ is the recovery rate of infected individuals,
		\item $\alpha > 0$ is the waning rate of the vaccine-induced immunity,
		\item $\gamma > 0$ is the waning rate of the infection-induced immunity.
	\end{itemize}
Within this SIRV model, the spread of the epidemic is described by the following system of differential equations:
	
\begin{eqnarray} \label{ode1}
	&& dS(t) = \left( \lambda - (\kappa + \rho) S(t) - \frac{\beta I(t)}{N(t)} S(t) + \alpha V(t) + \gamma R(t) \right) \, dt \nonumber \\
	&& dI(t) = \left( \frac{\beta I(t)}{N(t)} S(t) + \delta \frac{\beta I(t)}{N(t)} V(t) - (\kappa_1 + \theta) I(t) \right) \, dt \nonumber \\
	&& dR(t) = \left( \theta I(t) - (\kappa + \gamma) R(t) \right) \, dt \nonumber \\
	&& dV(t) = \left( \rho S(t) - (\kappa + \alpha) V(t) - \delta \frac{\beta I(t)}{N(t)} V(t) \right) \, dt,
\end{eqnarray}
and $dN(t) = \left( \lambda - \kappa N(t) - \left( \kappa_{1} - \kappa \right) I(t) \right) \, dt$ is the differential equation describing the population size over time.

\begin{remark} \label{rem-contact-rate}
	The incidence in system \eqref{ode1}, i.e. the number of new infections per unit of time, has the standard form $\beta I S / N$, and the parameter $\beta$ is understood as in \cite{Hethcote2000} -- $\beta$ is the average number of adequate contacts of one individual per unit of time, where a contact is called adequate if it is sufficient for the transmission of the disease. Consequently, $\beta I / N$ is the average number of adequate contacts with infective individuals per unit of time, i.e. the force of infection, while $(\beta I / N) S$ is the number of new infection cases per unit of time. In accordance with \cite{Hethcote2000}, the real parameter $\beta$ is therefore called the contact rate, which is in Subsection \ref{SIRV-time-changed} replaced by the stochastic process $\widetilde{\beta}$.
	
	The alternative mass-action incidence $\eta I S$ corresponds to $\beta = \eta N$, i.e. to a contact rate growing linearly with the population size. The data reviewed in \cite{Hethcote2000} indicate that the contact rate is nearly independent of the time-varying population size.
\end{remark}

\begin{remark} \label{rem-Ishikawa-deterministic}
	In particular, the system \eqref{ode1} differs from the deterministic SIRV model studied in \cite{im} in three respects. First, in \cite{im} $S(t)$, $I(t)$, $R(t)$ and $V(t)$ are observed as proportions of a population of constant size, so that $S(t) + I(t) + R(t) + V(t) = 1$ for all $t \geq 0$. In model \eqref{ode1} $S(t)$, $I(t)$, $R(t)$ and $V(t)$ are treated as numbers of individuals at time $t \geq 0$, and since the birth rate $\lambda$, the natural death rate $\kappa$ and the disease-induced death rate $\kappa_{1} \geq \kappa$ are three distinct quantities, the population size $N(t)$ varies over time. Second, the incidence in model \eqref{ode1} is of standard type, $\beta I S / N$, rather than of mass-action type $\beta S I$ as in \cite{im}. Third, the vaccination rate, which in \cite{im} is the control variable of an optimal vaccination problem, is in model \eqref{ode1} taken to be constant, since in this paper no optimization of the vaccination strategy is observed. The remaining parameters correspond one to one: the recovery rate $\theta$, the rates $\gamma$ and $\alpha$ of loss of the infection-induced and of the vaccine-induced immunity and the vaccine efficiency $\delta$ play the roles of $\gamma$, $\nu$, $\phi$ and $\sigma$ from \cite{im}, respectively.
	
	Furthermore, \cite{im} introduces stochasticity into the SIRV system through an additive Brownian perturbation of the contact rate $\beta$, whereas in the model developed in this paper stochasticity enters in a different way, as discussed in Remark \ref{rem-Ishikawa-stochastic}.
\end{remark}

Regarding long-term behaviour of the epidemic (extinction or persistence), the standard threshold quantity in compartmental epidemic models is the basic reproduction number $R_{0}$. It is defined as the expected number of secondary infections produced by a single infected individual introduced into a population at the disease-free state. The disease-free state is locally asymptotically stable if $R_{0} < 1$ (extinction scenario) and unstable if $R_{0} > 1$ (persistence scenario). For deterministic model \eqref{ode1} we compute $R_{0}$ by the next-generation method (see e.g. \cite{DHM1990,  DHR2010, vdDW2002}).

At the disease-free state $(S^{\ast}, I^{\ast}, R^{\ast}, V^{\ast})$ there are no infected, and consequently no recovered individuals $(I^{\ast} = R^{\ast} = 0)$. Therefore, the first and the fourth equation in \eqref{ode1} reduce to the linear ODEs
\begin{equation*}
	\frac{dS(t)}{dt} = \lambda - (\kappa + \rho) S(t) + \alpha V(t), \qquad \frac{dV(t)}{dt} = \rho S(t) - (\kappa + \alpha) V(t),
\end{equation*}
independent of the contact rate $\beta$. The unique equilibrium is
\begin{equation*} \label{DFE}
	E_{0} := \left( S^{\ast}, I^{\ast}, R^{\ast}, V^{\ast} \right) = \frac{\lambda}{\kappa} \left( \frac{\kappa + \alpha}{\kappa + \alpha + \rho}, \ 0, \ 0, \ \frac{\rho}{\kappa + \alpha + \rho} \right), \qquad N^{\ast} = S^{\ast} + I^{\ast} + R^{\ast} + V^{\ast} = \frac{\lambda}{\kappa},
\end{equation*}
and it is locally asymptotically stable since the Jacobian of the disease-free system has the negative eigenvalues $-\kappa$ and $-(\kappa + \rho + \alpha)$, and positive determinant $\kappa (\rho + \kappa + \alpha)$. The only infected compartment is $I$, with the rate of appearance of new infections $\mathcal{F}(S, I, R, V) = \beta \left( S + \delta V \right) I / N$, and the rate $\mathcal{V}(S, I, R, V) = \left( \kappa_{1} + \theta \right) I$ of transitions from $I$. The next-generation matrix $F V^{-1}$ (see \cite{vdDW2002}) is here of order one, with
\begin{equation} \label{R0-F-V}
F = \partial_{I} \mathcal{F} \big|_{E_{0}} = \beta \, \frac{S^{\ast} + \delta V^{\ast}}{N^{\ast}} = \beta \, \frac{\kappa + \alpha + \delta \rho}{\kappa + \alpha + \rho}, \quad V = \partial_{I} \mathcal{V} \big|_{E_{0}} = \kappa_{1} + \theta.
\end{equation} 
Therefore, the basic reproduction number of the model \eqref{ode1} is
\begin{equation} \label{R0-deterministic}
	R_{0} = \frac{\beta}{\kappa_{1} + \theta} \, \frac{\kappa + \alpha + \delta \rho}{\kappa + \alpha + \rho}.
\end{equation}
The first factor in \eqref{R0-deterministic} is the expected number of adequate contacts of an infected individual during its infectious period of expected length $1 / (\kappa_{1} + \theta)$. The second factor, coming from $F$ given in \eqref{R0-F-V}, is the fraction of the disease-free population effectively exposed to infection. It is equal to $1$ if $\rho = 0$ or $\delta = 1$, and it is an increasing function of $\delta$. Due to its dependence on parameters of the model, we observe \eqref{R0-deterministic} as a function $\beta \mapsto R_{0}(\beta)$.
	
For the long-term predictions, model \eqref{ode1} can be generalized by introducing the time-dependent contact rate $\widetilde{\beta}$ and by introducing the inflow into the susceptible class which, taking into account the birth and the death rates, ensures that at each time $t \geq 0$ the population size is bounded by the constant capacity $K > 0$, reflecting the limitations of the environmental resources. In this case, the reproduction number is $R_{0}(\widetilde{\beta})$ and it is itself time-dependent.

\subsection{Stochastic contact rate} \label{SIRV-time-changed}

To capture both the persistent evolution and the sudden bursts of the epidemic, the model for the contact rate $\widetilde{\beta}$ proposed here is driven by the time-changed L\'{e}vy noise $\mu$ introduced in Section \ref{subsec-preliminaries}:
\begin{equation} \label{SDE-beta}
	\widetilde{\beta}(t) = \widetilde{\beta}_{max} f\left( \beta(t) \right), \quad d\beta(t) = \beta(t) \, dt + \int\limits_{\mathbb{R}} \sigma_{t}(z, \beta(t)) \, \mu(dt, dz).
\end{equation}
The driving noise $\mu$ is addressing both continuous fluctuations and sporadic abrupt jumps in contact behavior, governed by $\Lambda^{B}$ and $\Lambda^{H}$, respectively. By construction, the underlying (latent) process $\beta$ given in \eqref{SDE-beta} is a semimartingale (it is the sum of an absolutely continuous drift and a stochastic integral with respect to the martingale random measure $\mu$). In definition of the stochastic contact rate $\widetilde{\beta}$, factor $\widetilde{\beta}_{max} > 0$ is the ceiling of the contact rate, i.e. the largest average number of adequate contacts per individual per unit of time (a constant of the same nature as the parameter $\beta$ from Subsection \ref{SIRV-deterministic}), and $f \colon \mathbb{R} \to [0, 1]$ is a twice continuously differentiable function assumed to satisfy the following requirements:
\begin{itemize}
	\item[(f1)] $0 < f_{min} \leq f(x) \leq 1$, for all $x \in \mathbb{R}$,
	\item[(f2)] $\sup\limits_{x \in \mathbb{R}} |f^{\prime}(x)| \leq K_{f^{\prime}}$, $K_{f^{\prime}} > 0$,
	\item[(f3)] $\sup\limits_{x \in \mathbb{R}} |f^{\prime \prime}(x)| \leq K_{f^{\prime \prime}}$, $K_{f^{\prime \prime}} > 0$.
\end{itemize}
Assumptions \it{(f1)-(f3)} ensure that the contact rate process $\widetilde{\beta} = (\widetilde{\beta}(t), \, t \geq 0)$, defined on the complete filtered probability space $\left( \Omega, \mathcal{F}, P, \mathbb{F} \right)$ introduced in Subsection \ref{subsec-preliminaries}), is bounded from below by $f_{min} \widetilde{\beta}_{max}$ and from above by $\widetilde{\beta}_{max}$, thus providing the (rough) estimates used in the proofs in the following Sections. 

The motivation for studying this particular contact rate model is to explore how the time-scale disturbances in the noise $\mu$ influence the spread of the disease. To the best of our knowledge, this particular time-change method is novel in modeling of epidemics.

Since the contact rate process $\widetilde{\beta}$ is just a functional transformation of the process $\beta$, its well-posedness follows from that of process $\beta$. According to \cite{GS2} (Section 6, p. 1699), the conditions for existence and uniqueness of a strong solution of the SDE \eqref{SDE-beta}, written in integral form as
\begin{equation*} \label{SDE-beta-integral}
	\beta(t) = \beta(0) + \int\limits_{0}^{t} \beta(s) \, ds + \int\limits_{0}^{t}  \sigma_{s}(0, \beta(s)) \, dB(s) + \int\limits_{0}^{t} \int\limits_{\mathbb{R}_{0}} \sigma_{s}(z, \beta(s)) \, \widetilde{H}(ds, dz),
\end{equation*}
rest on four requirements. The first two concern the initial value and the regularity of the jump coefficient in the state variable:
\begin{itemize}
	\item[($e1$)] $\beta(0)$ is $\mathcal{F}_{0}$-measurable and square integrable,
	\item[($e2$)] the jump coefficient $\sigma_{t}(z,b)$ is Lipschitz continuous, i.e. there exists a constant $k_{1} > 0$ such that for all $t \geq 0$ and $b_{1}, b_{2} \in \mathbb{R}$
		\begin{align*}
			& \left|\sigma_{t}(0,b_{1})-\sigma_{t}(0,b_{2})\right| \le k_{1} \left|b_{1}-b_{2}\right|, & \\
			& \left|\sigma_{t}(z,b_{1})-\sigma_{t}(z,b_{2})\right| \le k_{1} \left|z(b_{1}-b_{2})\right|, \quad z \neq 0. &
		\end{align*}
		\end{itemize}
Quadratic variation and the predictable quadratic variation for the process $\beta$ are (see Appendix \ref{app:noise}):
\begin{align} \label{beta-qv}
	[\beta, \beta]_{t} & = \int_{0}^{t} \sigma_{s}^{2}(0, \beta(s-)) \lambda_{s}^{B} \, ds + \sum\limits_{0 < s \leq t} \left( \Delta \beta(s) \right)^{2}, \nonumber \\
	\langle \beta, \beta \rangle_{t} & = \int_{0}^{t} \sigma_{s}^{2}(0, \beta(s-)) \lambda_{s}^{B} \, ds + \int_{0}^{t} \int_{\mathbb{R}_{0}} \sigma_{s}^{2}(z, \beta(s-)) \lambda_{s}^{H} \nu(dz) \, ds,
\end{align}
with $\Delta \beta(s) = \int_{\mathbb{R}_{0}} \sigma_{s}(z, \beta(s-)) \, H(\{s\}, dz)$, both finite $P$-a.s. for every $t > 0$.

The remaining two requirements for existence and uniqueness of solution of the SDE \eqref{SDE-beta} address the integrability of the jump coefficient with respect to random measure $\mu$ (otherwise the stochastic integral in \eqref{SDE-beta} is not well defined):
\begin{itemize}
	\item[$\it{(\sigma1)}$] $\sigma_{t}(z, \beta(t))$ is $\mathbb{F}$-predictable,
	\item[$\it{(\sigma2)}$] there exist a deterministic measurable function $g \colon [0, \infty \rangle \times \mathbb{R} \to [0, \infty \rangle$ which dominates the jump coefficient uniformly in the state variable, i.e.
	\begin{equation} \label{sigma-envelope}
		\left| \sigma_{t}(z, b) \right| \leq g_{t}(z), \qquad \text{for all } t \geq 0, \ z \in \mathbb{R}, \ b \in \mathbb{R}, \quad P\text{-a.s.},
	\end{equation}
	and which is locally bounded in the sense that for every $t > 0$
	\begin{equation} \label{sigma-envelope-moments}
		g_{t}^{\ast}(0) := \sup\limits_{s \leq t} g_{s}(0) < \infty,
		\qquad
		G_{2}^{\ast}(t) := \sup\limits_{s \leq t} \int\limits_{\mathbb{R}_{0}} g_{s}^{2}(z) \, \nu(dz) < \infty.
	\end{equation}
	In epidemiological interpretation, $g_{t}(z)$ describes the largest response of the process $\beta$ to a shock of size $z$ occurring at time $t$.
\end{itemize}
	
Therefore, under requirements $\it{(f1)-(f3), (e1)-(e2)}$ and $\it{(\sigma1)-(\sigma2)}$, the process $\widetilde{\beta}$ is well defined.

\begin{remark} \label{rem-affine-drift}
	The identity drift in the SDE \eqref{SDE-beta} is chosen for simplicity. More generally, the drift can be affine, i.e. of the form $a + b \beta(t)$ with $a, b \in \mathbb{R}$, which for $b < 0$ includes the mean-reverting case. The arguments used to prove the results in this paper do not depend on the particular choice of the affine drift -- the estimates in Sections \ref{sec:EUPsolution}, \ref{sec:extinction} and \ref{sec:persistence} use only the boundedness of the function $f$ and of its derivatives, and the drift of the process $\beta$ appears only in application of the It\^{o} formula in the term $\beta \cdot f^{\prime}(\beta)$ which, in case of the general affine drift, becomes $f^{\prime}(\beta) (a + b \beta)$, and is controlled by assumption $\it{(f2)}$ and additional assumption $\it{(f4)}$ introduced in Section \ref{sec:extinction}. The identity drift is thus the simplest representative of the affine class, and it keeps the constants in the thresholds in Sections \ref{sec:extinction} and \ref{sec:persistence} more tractable.
\end{remark}

\begin{exa}[Function $f$ and jump coefficient $\sigma$] \label{ex-f-sigma}
\begin{itemize}
	\item[]
	\item[(a)] An example of the function satisfying $\it{(f1)-(f3)}$, and also the condition $\it{(f4)}$ introduced in Section \ref{sec:extinction}, is 
	\begin{equation} \label{f-example}
		f(x) = f_{min} + \left( 1 - f_{min} \right) \frac{1 + \tanh (x)}{2}, \qquad x \in \mathbb{R}, \quad 0 < f_{min} < 1.
	\end{equation}
	Its derivatives are $f^{\prime}(x) = \frac{1 - f_{min}}{2} \operatorname{sech}^{2} (x)$ and $f^{\prime \prime}(x) = - \left( 1 - f_{min} \right) \operatorname{sech}^{2} (x) \cdot \tanh (x)$, yielding
	\begin{equation*}
		K_{f^{\prime}} = \frac{1 - f_{min}}{2}, \qquad K_{f^{\prime \prime}} = \frac{2}{3 \sqrt{3}} \left( 1 - f_{min} \right), \qquad K_{x f^{\prime}} = \frac{1 - f_{min}}{2} \max_{x \in \mathbb{R}} \left| x \operatorname{sech}^{2} x \right|.
	\end{equation*}
	Note that $f$ approaches its bounds $f_{min}$ and $1$ exponentially fast as $x \to -\infty$ and $x \to \infty$, respectively, without attaining them, so that $\widetilde{\beta}_{max} f_{min} < \widetilde{\beta}(t) < \widetilde{\beta}_{max}$. 
	Another example of a function satisfying $\it{(f1)-(f4)}$ is
	\begin{equation*} \label{f-example-bell}
		f(x) = f_{min} + \left( 1 - f_{min} \right) e^{-x^{2}}, \qquad x \in \mathbb{R}, \quad 0 < f_{min} < 1,
	\end{equation*}
	with derivatives $f^{\prime}(x) = - 2 \left( 1 - f_{min} \right) x e^{-x^{2}}$ and $f^{\prime \prime}(x) = 2 \left( 1 - f_{min} \right) \left( 2 x^{2} - 1 \right) e^{-x^{2}}$. Therefore,
	\begin{equation*}
		K_{f^{\prime}} = \sqrt{\frac{2}{e}} \left( 1 - f_{min} \right), \qquad K_{f^{\prime \prime}} = 2 \left( 1 - f_{min} \right), \qquad K_{x f^{\prime}} = \frac{2}{e} \left( 1 - f_{min} \right).
	\end{equation*}
	All the requirements on function $f$ are satisfied, but this choice is not a reasonable one from epidemiological point of view. Namely, $f$ is not monotone -- it achieves its maximal value $1$ at $x = 0$ and decreases as $\left| x \right| \to \infty$, so the contact rate is at its peak when $\beta(t) = 0$ and decays towards $\widetilde{\beta}_{max} f_{min}$ as the process $\beta$ moves away from $0$, in either direction. Therefore, under this choice for $f$, the most cautious and the most reckless behavioral scenarios in population would produce similar number of contacts, which is not the case in reality.
	 
	\item[(b)] For the jump coefficient, a simple state-independent choice satisfying $\it{(e2),(\sigma1)}$ and $\it{(\sigma2)}$ is
	\begin{equation} \label{sigma-example-affine}
		\sigma_{t}(z, b) = s_{0} \, \mathbf{1}_{\{ z = 0 \}} + s_{1} z \, \mathbf{1}_{\{ z \neq 0 \}}, \qquad s_{0}, s_{1} \in \mathbb{R},
	\end{equation}
	with $g_{t}(z) = |s_{0}| \mathbf{1}_{\{ z = 0 \}} + |s_{1}| |z| \mathbf{1}_{\{ z \neq 0 \}}$, for which $G_{2} = s_{1}^{2} \int_{\mathbb{R}_{0}} z^{2} \nu(dz)$ and $G_{4} = s_{1}^{4} \int_{\mathbb{R}_{0}} z^{4} \nu(dz)$. Therefore, the additional condition $\it{(\sigma2^{\ast})}$ introduced in Section \ref{sec:extinction} holds as soon as $\nu$ has a finite fourth moment. A state-dependent example is $\sigma_{t}(z, b) = \left( s_{0} \, \mathbf{1}_{\{ z = 0 \}} + s_{1} z \, \mathbf{1}_{\{ z \neq 0 \}} \right) h(b)$, where the function $h$ is bounded and Lipschitz continuous.
\end{itemize}
\end{exa}
\begin{figure}[H]
	\centering
	\includegraphics[width=\textwidth]{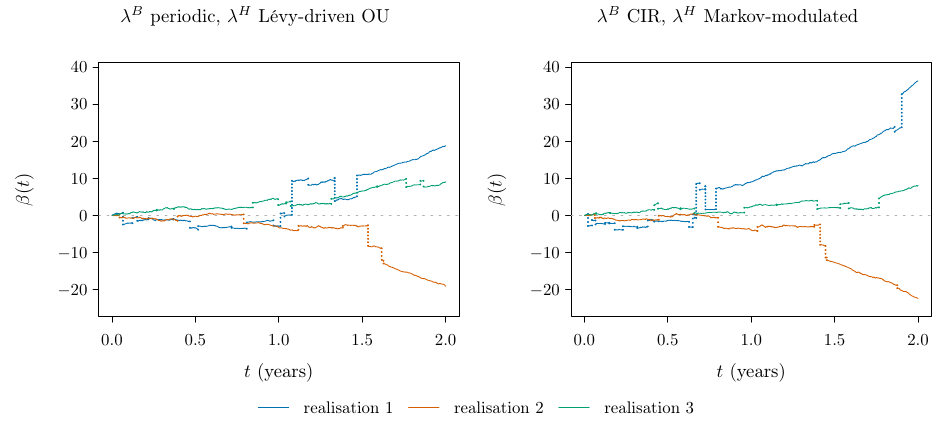}
	\caption{Sample paths of the underlying process $\beta$ with the jump coefficient \eqref{sigma-example-affine}, for three realisations of the noise and the two stochastic clocks with the intensities from Figure~\ref{fig:clocks-intensity}. Parameter values are given in Section~\ref{sec:simulation}.} \label{fig:beta}
\end{figure}

\begin{figure}[htbp]
	\centering
	\includegraphics[width=\textwidth]{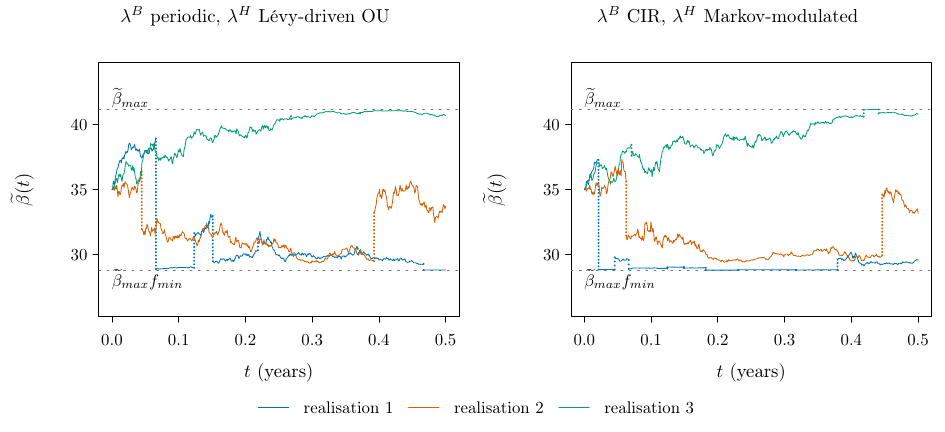}
	\caption{Sample paths of the contact rate $\widetilde{\beta} = \widetilde{\beta}_{max} f(\beta)$, with $f$ given by \eqref{f-example}, corresponding to the paths of $\beta$ from Figure~\ref{fig:beta}, in the extinction regime ($\overline{\mathcal{R}}_{0} = 0.7$, see Section~\ref{sec:extinction}) over the first half year. Parameter values are given in Section~\ref{sec:simulation}.} \label{fig:btilde-extinction}
\end{figure}
\begin{figure}[htbp]
	\centering
	\includegraphics[width=\textwidth]{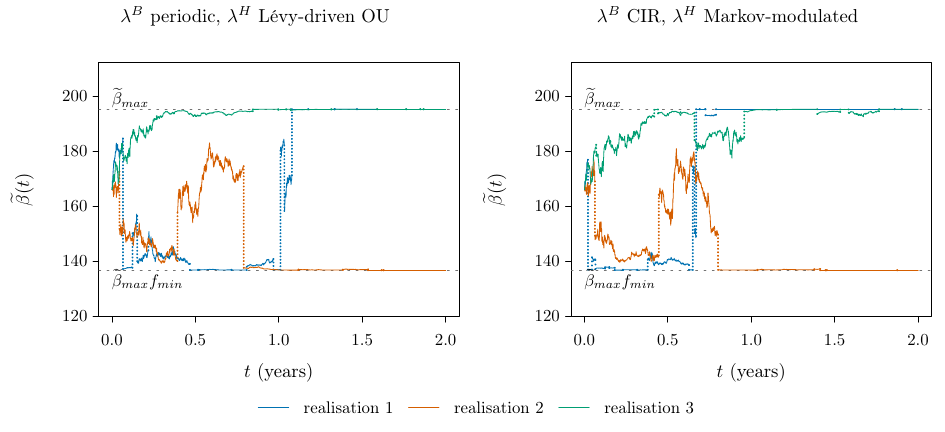}
	\caption{Sample paths of the contact rate $\widetilde{\beta} = \widetilde{\beta}_{max} f(\beta)$, with $f$ given by \eqref{f-example}, corresponding to the paths of $\beta$ from Figure~\ref{fig:beta}, in the persistence regime ($\underline{\mathcal{R}}_{0} = 1.16$ and $\overline{\mathcal{R}}_{0} = 3.32$, see Section~\ref{sec:persistence}) over two years. Parameter values are given in Section~\ref{sec:simulation}.} \label{fig:btilde-persistence}
\end{figure}
\begin{figure}[htbp]
	\centering
	\includegraphics[width=\textwidth]{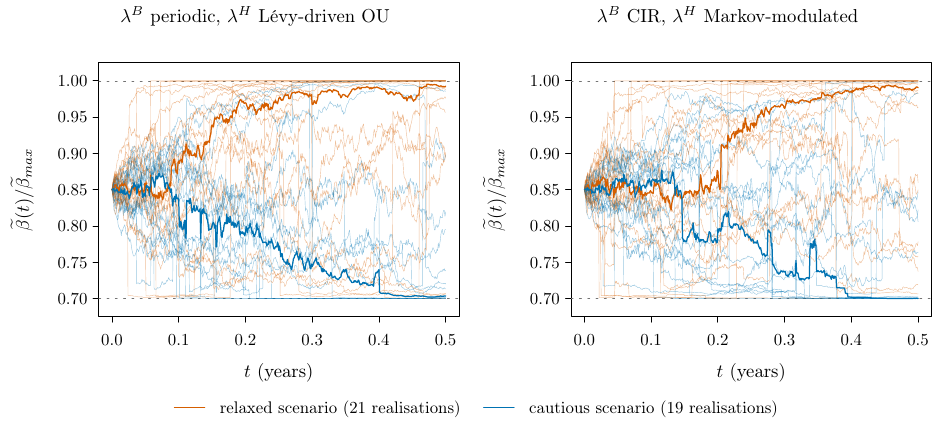}
	\caption{The relative contact rate $\widetilde{\beta}(t) / \widetilde{\beta}_{max} = f(\beta(t))$ for $f$ given by \eqref{f-example}: $40$ realisations of the noise and the two stochastic clocks (left and right). The realisations are coloured according to the level at which the contact rate eventually settles, $\widetilde{\beta}_{max}$ (relaxed scenario) or $\widetilde{\beta}_{max} f_{min}$ (cautious scenario), and the thick lines are the pointwise medians of the two groups. Parameter values are given in Section~\ref{sec:simulation}.} \label{fig:ensemble-beta}
\end{figure}

\subsection{Random SIRV epidemic model} \label{subsec-stochastic-SIRV}

In the deterministic model \eqref{ode1}, $S$, $I$, $R$ and $V$ are deterministic functions of time. Replacing the constant contact rate $\beta$ by the stochastic process $\widetilde{\beta}$ defined in \eqref{SDE-beta} turns them into stochastic processes, and yields a more realistic description of the epidemic dynamics.

Note that in reality the population size is not constant, but rather follows a dynamic described by the process $N = \left( N(t), \, t \geq 0 \right)$. Also, it is natural to assume that at each $t \geq 0$ the population size $N(t)$ is bounded by the environmental carrying capacity $K$. Moreover, $N(t)$ is influenced by the birth rate $\lambda > 0$, the natural death rate $\kappa > 0$ and the disease-induced death-rate $\kappa_{1} \geq \kappa$. To embed these influences into the model, we introduce the constant
\begin{equation} \label{H-population-size}
	H := \frac{\kappa}{\lambda} K, \qquad K > 0,
\end{equation}
and replace the constant birth inflow $\lambda$ in system \eqref{ode1} by $\lambda H = \kappa K$, meaning that susceptibles are added at exactly the rate at which natural deaths remove individuals when the population is at its capacity $K$. If $\kappa < \lambda$, $H$ is a fraction of the environmental capacity $K$. On the other hand, $\kappa > \lambda$ would require more individuals than the environment can hold, and therefore is not considered in this modeling framework.

By introducing the population dynamics described above into the model \eqref{ode1}, we obtain the compartmental model with the random coefficient $\widetilde{\beta}$, coupled with the SDE \eqref{SDE-beta} for the contact rate, which we call {\it the random SIRV model}:
\begin{eqnarray}
	&& dS(t) = \left(\lambda H - (\rho + \kappa) S(t) - \frac{\widetilde{\beta}(t) I(t)}{N(t)} S(t) + \alpha V(t) + \gamma R(t) \right) \, dt \nonumber \\
	&& dI(t) = \left( \frac{\widetilde{\beta}(t) I(t)}{N(t)} S(t) + \delta \frac{\widetilde{\beta}(t) I(t)}{N(t)} V(t) - (\kappa_1 + \theta) I(t) \right) \, dt \nonumber \\
	&& dR(t) = \left( \theta I(t) - (\kappa + \gamma) R(t) \right) \, dt \label{sde1} \\
	&& dV(t) = \left( \rho S(t) - (\kappa + \alpha) V(t) - \delta \frac{\widetilde{\beta}(t) I(t)}{N(t)} V(t) \right) \, dt \nonumber \\
	&& d\widetilde{\beta}(t) = \widetilde{\beta}_{max} \, d f\left( \beta(t) \right), \nonumber
\end{eqnarray}
where, by the It\^{o} formula for semimartingales applied to $f(\beta)$, we obtain the SDE for its functional transformation $f(\beta)$: 
\begin{align} \label{dfbeta}
	d f\left( \beta(t) \right) & = \left( f^{\prime}\left( \beta(t-) \right) \beta(t)
		+ \frac{1}{2} f^{\prime \prime}\left( \beta(t-) \right) \sigma_{t}^{2}\left( 0, \beta(t-) \right) \lambda_{t}^{B} \right. \\
	& \left. + \lambda_{t}^{H} \int\limits_{\mathbb{R}_{0}} \left( f\left( \beta(t-) + \sigma_{t}(z, \beta(t-)) \right)
	- f\left( \beta(t-) \right) - f^{\prime}\left( \beta(t-) \right) \sigma_{t}(z, \beta(t-)) \right) \nu(dz) \right) dt \nonumber \\
	& + \int\limits_{\mathbb{R}} \left[ \mathbf{1}_{\{ z = 0 \}} \, f^{\prime}\left( \beta(t-) \right) \sigma_{t}\left( 0, \beta(t-) \right)
	+ \mathbf{1}_{\{ z \neq 0 \}} \left( f\left( \beta(t-) + \sigma_{t}(z, \beta(t-)) \right) - f\left( \beta(t-) \right) \right) \right] \mu(dt, dz). \nonumber
\end{align}
Note that the fifth equation of system \eqref{sde1} is the transformed version of the SDE \eqref{SDE-beta} for the process $\beta$, and its solution enters the SIRV system as a random, c\`{a}dl\`{a}g and bounded coefficient. This is the sense in which the model is random, rather than stochastic -- given the path of $\widetilde{\beta}$, the compartments evolve deterministically, and the randomness of $S$, $I$, $R$ and $V$ is inherited from $\widetilde{\beta}$. All the parameters of the model are interpreted in Subsection \ref{SIRV-deterministic}.

\begin{remark} \label{rem-Ishikawa-stochastic}
	In \eqref{sde1} randomness enters the model only through the stochastic contact rate $\widetilde{\beta}$. This differs from the approach used in \cite{im}, whose deterministic SIRV model inspired \eqref{ode1}, where the contact rate is perturbed by a Gaussian white noise $\beta \, dt \mapsto \beta \, dt + \varepsilon \, dW(t)$, where $W$ is the standard Brownian motion. Furthermore, e.g. in \cite{Tymoshenko-2025}, each compartment equation is driven by its own Brownian motion, while in \cite{Tymoshenko-2026} a self-exciting jump noise of Hawkes type is added to each compartment. In random system \eqref{sde1} the contact rate is a stochastic process driven by the time-changed L\'{e}vy noise, and it enters the compartmental ODEs as a random coefficient. Therefore, there is no choice of $f$, $\sigma$ and $\Lambda$ that turns \eqref{sde1} into the model studied in \cite{im}.
\end{remark}

Note that by summing the first four equations of \eqref{sde1} and using relation $\lambda H = \kappa K$ from \eqref{H-population-size}, the terms containing $\widetilde{\beta}$ cancel, and the population size $N(t) = S(t) + I(t) + R(t) + V(t)$ satisfies the following differential equation:
\begin{equation} \label{N-equation}
	\frac{d N(t)}{dt} = \lambda H - \kappa \left( S(t) + R(t) + V(t) \right) - \kappa_{1} I(t) = \kappa \left( K - N(t) \right) - \left( \kappa_{1} - \kappa \right) I(t).
\end{equation}
The rational assumption is that $0 < N(0) < K$, where $K$ is the environmental capacity of the population.

\begin{remark} \label{rem-R0}
	The basic reproduction number \eqref{R0-deterministic} does not depend on the inflow $\lambda$ into the susceptible class, and the disease-free state does not depend on the contact rate. Replacing $\lambda$ by $\kappa K$ therefore leaves the function $\beta \mapsto R_{0}(\beta)$ unchanged, and the basic reproduction number of the random model \eqref{sde1} is the stochastic process obtained by composing this function with the contact rate $\widetilde{\beta}$, i.e.
	\begin{equation*} \label{R0-definition}
		\mathcal{R}_{0}(t) := R_{0}\left( \widetilde{\beta}(t) \right) = \frac{\widetilde{\beta}(t)}{\kappa_{1} + \theta} \,  \, \frac{\kappa + \alpha + \delta \rho}{\kappa + \alpha + \rho} = \frac{\widetilde{\beta}_{max} f\left( \beta(t) \right)}{\kappa_{1} + \theta} \,  \, \frac{\kappa + \alpha + \delta \rho}{\kappa + \alpha + \rho}, \qquad t \geq 0 .
	\end{equation*}
	Since it is a linear transformation of the contact rate $\widetilde{\beta}$, $\mathcal{R}_{0}$ is a c\`{a}dl\`{a}g, $\mathbb{F}$-adapted process. Its interpretation at time $t$ is as follows: $\mathcal{R}_{0}(t)$ is the expected number of secondary infections that a single infected individual, introduced at time $t$ into the disease-free population, would produce if the contact rate is held at its current level $\widetilde{\beta}(t)$. 
	
	Since by $\it{(f1)}$ the function $f$ takes values in $[f_{min}, 1]$ and since $\frac{\kappa + \alpha + \delta \rho}{\kappa + \alpha + \rho} \in [\delta, 1]$, the process $\mathcal{R}_{0}$ is bounded:
	\begin{equation*} \label{R0-envelopes}
		\underline{\mathcal{R}}_{0} := \frac{\delta \widetilde{\beta}_{max} f_{min}}{\kappa_{1} + \theta} \ \leq \ \mathcal{R}_{0}(t) \ \leq \ \frac{\widetilde{\beta}_{max}}{\kappa_{1} + \theta} =: \overline{\mathcal{R}}_{0}.
	\end{equation*}
\end{remark}

\begin{remark} \label{rem5.1}
	Beside the contact rate $\widetilde{\beta}$ , other parameter of the model \eqref{sde1} also could be treated as a bounded deterministic or stochastic functions of time. The analytical framework remains similar to the time-independent case, but the transition rates between compartments would be time-dependent and the functional form of each parameter should reflect its biological meaning, supported by data or by modeling assumptions. 
	
	Furthermore, the randomness of the system \eqref{sde1} can be enhanced by introducing noises in the equations for compartments $S$, $I$, $R$ and $V$ as well. These noises could be mutually independent or correlated, as well as independent or correlated with the contact rate $\widetilde{\beta}$. This turns the random model into a proper stochastic model. The analysis of such systems is more complex, although the techniques remain the same (see Remark \ref{rem4.1}). This approach is left for future research -- the emphasis of the present paper is on studying the effect of the absolutely continuous stochastic clock $\Lambda$ on the extinction and persistence scenarios within model \eqref{sde1}.
\end{remark}

\section{Existence and uniqueness of the global solution} \label{sec:EUPsolution}

In this section we prove existence and uniqueness of the positive global solution of the random SIRV model \eqref{sde1}, and show that the population size stays within the environmental capacity. It is convenient to write the system \eqref{sde1} in the following compact form:
\begin{equation} \label{sde-vector-form}
	X(t) = X(0) + \int\limits_{0}^{t} a(s, X(s)) \, ds +  \int\limits_{0}^{t} \int\limits_{\mathbb{R}} b(s, z) \, \mu(ds, dz).
\end{equation}
In this form, $X(t) = \left( S(t), I(t), R(t), V(t), \widetilde{\beta}(t) \right)$, and the drift $a \colon \mathbb{R}_{+} \times \mathcal{D} \to \mathbb{R}^{5}$, where $\mathbb{R}_{+}$ is the set of all non-negative real numbers and $\mathcal{D} := \left\{ \mathbf{x} \in \mathbb{R}^{5} : \, x_{1} + x_{2} + x_{3} + x_{4} \neq 0 \right\}$, is measurable in the time variable and continuously differentiable in the space variable. As usual, the dependence on $\omega \in \Omega$ is suppressed in the notation, and the drift is explicitly given by
\begin{equation} \label{sde1-vector-form-a}
	a(t, X(t))= \begin{bmatrix}
		\lambda H - (\rho + \kappa) S(t) - \displaystyle\frac{\widetilde{\beta}(t) I(t)}{N(t)} S(t) + \alpha V(t) +\gamma R(t) \\ \\
		\displaystyle\frac{\widetilde{\beta}(t) I(t)}{N(t)} S(t) + \delta \displaystyle\frac{\widetilde{\beta}(t) I(t)}{N(t)} V(t) - (\kappa_1 + \theta)I(t) \\ \\
		\theta I(t) - (\kappa + \gamma) R(t) \\ \\
		\rho S(t) - (\kappa + \alpha) V(t) - \delta \displaystyle\frac{\widetilde{\beta}(t) I(t)}{N(t)} V(t) \\ \\
		a_{5}(t)
	\end{bmatrix}.
\end{equation}
Here, $a_{i}(t, X(t))$, $i \in \{1, \ldots, 5\}$, denotes the $i$-th component of the function $a$. The first four components are deterministic functions of the state of the process $X$, and they depend on time and on $\omega$ only through the fifth coordinate $\widetilde{\beta}(t)$. The fifth component, on the other hand, does not depend on the state of the process $X$, i.e. $a_{5}(t, X(t)) = a_{5}(t)$. Like the jump coefficient $b$, it is determined by the underlying process $\beta$, which solves \eqref{SDE-beta} independently of the compartments. Explicitly, the fifth component of the drift $a$ is given by
\begin{align} \label{a5-definition}
	a_{5}(t) & = \widetilde{\beta}_{max} \left( f^{\prime}\left( \beta(t-) \right) \beta(t)
	+ \frac{1}{2} f^{\prime \prime}\left( \beta(t-) \right) \sigma_{t}^{2}\left( 0, \beta(t-) \right) \lambda_{t}^{B} \right. \nonumber \\
	& \qquad \left. + \lambda_{t}^{H} \int\limits_{\mathbb{R}_{0}} \left( f\left( \beta(t-) + \sigma_{t}(z, \beta(t-)) \right) \nonumber
	- f\left( \beta(t-) \right) - f^{\prime}\left( \beta(t-) \right) \sigma_{t}(z, \beta(t-)) \right) \nu(dz) \right). \nonumber
\end{align}
Further, the jump coefficient $b$, with components $b_{i}$, $i \in \{1, \ldots, 5\}$, is $\mathbb{F}$-predictable and its only non-zero component $b_{5}$ is given by
\begin{equation*} \label{b5-definition}
	b_{5}(t, z) = \widetilde{\beta}_{max} \left( \mathbf{1}_{\{ z = 0 \}} \, f^{\prime}\left( \beta(t-) \right) \sigma_{t}\left( 0, \beta(t-) \right)
	+ \mathbf{1}_{\{ z \neq 0 \}} \left( f\left( \beta(t-) + \sigma_{t}(z, \beta(t-)) \right) - f\left( \beta(t-) \right) \right) \right).
\end{equation*}

Before stating the theorem on the existence and the uniqueness of the positive global solution of system \eqref{sde-vector-form}, let us recall the framework in which the main result of this section, and all the results that follow, are obtained. We work on the complete filtered probability space $\left( \Omega, \mathcal{F}, P, \mathbb{F} \right)$ from Subsection \ref{subsec-preliminaries}, with the right-continuous filtration $\mathbb{F} = \left( \mathcal{F}_{t}, \, t \geq 0 \right)$ generated by the random measure $\mu$. Further, it is assumed that $\it{(l1)-(l3), (f1)-(f3), (e1)-(e2)}$ and $\it{(\sigma1)-(\sigma2)}$ hold.

\begin{theorem} \label{thm-global-positive-solution}
	Let the initial values of all compartments of the random SIRV model \eqref{sde-vector-form} be $\mathcal{F}_{0}$-measurable and satisfy
	\begin{equation} \label{initial-condition}
		\left( S(0), I(0), R(0), V(0) \right) \in \langle 0, \infty \rangle^{4},
		\qquad 0 < N(0) \leq K, \qquad P\text{-a.s.}
	\end{equation}
	Then the system \eqref{sde-vector-form} has a unique global solution
	$X = ((S(t), I(t), R(t), V(t), \widetilde{\beta}(t)), \, t \geq 0)$, whose first four components are $\mathbb{F}$-adapted and absolutely continuous and whose fifth component is $\mathbb{F}$-adapted and c\`{a}dl\`{a}g, and $P$-a.s. for every $t \geq 0$:
	\begin{itemize}
		\item[(i)] $S(t), I(t), R(t), V(t) > 0$,
		\item[(ii)] the population size satisfies
		\begin{equation} \label{N-bounds}
			\frac{\kappa K}{\kappa_{1}} + \left( N(0) - \frac{\kappa K}{\kappa_{1}} \right) e^{-\kappa_{1} t}
			\ \leq \ N(t) \ \leq \ K + \left( N(0) - K \right) e^{-\kappa t} ,
		\end{equation}
		and consequently
		\begin{equation} \label{N-boundaries}
			\frac{\kappa}{\kappa_{1}} K \ \leq \ \liminf\limits_{t \to \infty} N(t) \ \leq \ \limsup\limits_{t \to \infty} N(t) \ \leq \ K .
		\end{equation}
	\end{itemize}
	In particular, the set
	\begin{equation} \label{gamma-set}
		\Gamma = \left\{( S, I, R, V, \widetilde{\beta}) \in \langle 0, K]^{4} \times \left[ \widetilde{\beta}_{max} f_{min}, \widetilde{\beta}_{max} \right] :
		\frac{\kappa}{\kappa_{1}} K \leq S + I + R + V \leq K \right\}
	\end{equation}
	is positively invariant for the process $X$, i.e. every solution started in $\Gamma$ $P$-a.s. stays in $\Gamma$ for all $t > 0$.
\end{theorem}

\begin{remark}
Epidemiologically, the theorem says that the model never produces a negative compartment and never lets the population exceed the environmental capacity $K$ or fall below $\kappa K / \kappa_{1}$. The lower bound $\kappa K / \kappa_{1}$ is the population size at which the inflow $\lambda H = \kappa K$ balances the deaths in the worst case, in which every individual is infected but the population does not die out, while the upper bound $K$ is the environmental carrying capacity. Note that the lower bound in \eqref{N-bounds} degenerates when there is no additional mortality $\kappa_{1} - \kappa$ induced by the disease. Namely, if $\kappa_{1} = \kappa$, the two bounds in \eqref{N-bounds} coincide, \eqref{N-equation} reduces to $\frac{d N}{dt} = \kappa (K - N)$, and $N(t) = K + (N(0) - K) e^{-\kappa t} \to K$ as $t \to \infty$.
\end{remark}

\begin{proof} \phantom{a} \newline \vspace*{-0.4cm}
	
	\noindent\textbf{The underlying process $\beta$ and the contact rate $\widetilde{\beta}$} \\
	The coefficients of the SDE \eqref{SDE-beta} satisfy the requirements given in \cite{GS2}, Section $6$, recalled in Subsection \ref{SIRV-time-changed}: the drift is Lipschitz, the initial value is square integrable by $\it{(e1)}$, and the jump coefficient is Lipschitz in the state variable by $\it{(e2)}$, $\mathbb{F}$-predictable by $\it{(\sigma1)}$ and square integrable with respect to the stochastic clock on bounded time-intervals -- by $\it{(\sigma2)}$ and $\it{(l1)}$, for every fixed $b \in \mathbb{R}$ and every $t > 0$,
	\begin{equation} \label{e3-implied}
		\int\limits_{0}^{t} \int\limits_{\mathbb{R}} \left| \sigma_{s}(z, b) \right|^{2} \Lambda(ds, dz)
		\leq \left( g_{t}^{\ast}(0) \right)^{2} \Lambda_{t}^{B} + G_{2}^{\ast}(t) \, \widehat{\Lambda}_{t}^{H} ,
	\end{equation}
	and the right-hand side is $P$-a.s. finite, with a finite expectation by $\it{(l3)}$. Since $\left| \sigma_{t}(z, b) \right| \leq g_{t}(z)$ holds uniformly in the state variable $b$ (see \eqref{sigma-envelope} from $\it{(\sigma2)}$), \eqref{e3-implied} remains valid with $b$ replaced by $\beta(s-)$. Hence, the SDE \eqref{SDE-beta} for the process $\beta$ has a unique strong solution on $[0, T]$ for every $T > 0$, and by uniqueness these solutions are restrictions of a single process $\beta$ on $[0, \infty \rangle$. 
	
	Being the sum of an absolutely continuous drift and of a stochastic integral with respect to the martingale random measure $\mu$, process $\beta$ is an $\mathbb{F}$-adapted c\`{a}dl\`{a}g semimartingale. It is bounded on compact time-intervals path by path,
	\begin{equation} \label{beta-bound}
		C_{T}^{(\beta)} := \sup\limits_{0 \leq t \leq T} \left| \beta(t) \right| < \infty, \quad P\text{-a.s., for every } T > 0,
	\end{equation}
	where $C_{T}^{(\beta)}$ doesn't have to be bounded uniformly in $\omega$.
	
	However, the Lyapunov argument used below in this proof requires a bound on $E[\sup_{t \leq T} \left| \beta(t) \right|]$ (see \eqref{L-expectation}). Since the SDE \eqref{SDE-beta} for process $\beta$ has the identity drift, the It\^{o} formula applied to $e^{-t} \beta(t)$ gives the explicit solution
	\begin{equation} \label{beta-explicit}
		\beta(t) = e^{t} \left( \beta(0) + Z(t) \right), \qquad Z(t) := \int\limits_{0}^{t} e^{-s} \int\limits_{\mathbb{R}} \sigma_{s} \left( z, \beta(s-) \right) \mu(ds, dz), \qquad t \geq 0 .
	\end{equation}
	Since $e^{-2s} \leq 1$, the process $Z$ is a square-integrable martingale on every bounded time-interval and, by the isometry \eqref{mu-isometry} from Appendix \ref{app:noise},
	\begin{equation*}
		E \left[ Z^{2}(T) \right] = E \left[ \int\limits_{0}^{T} \int\limits_{\mathbb{R}} e^{-2s} \sigma_{s}^{2} \left( z, \beta(s-) \right) \Lambda(ds, dz) \right]
		\leq \left( g_{T}^{\ast}(0) \right)^{2} E \left[ \Lambda_{T}^{B} \right] + G_{2}^{\ast}(T) \, E \left[ \widehat{\Lambda}_{T}^{H} \right] < \infty .
	\end{equation*}
	Therefore, by $(x + y)^{2} \leq 2 \left( x^{2} + y^{2} \right)$, Doob's $L^{2}$ inequality and $(e1)$, for every $T > 0$
	\begin{equation} \label{beta-moment}
		E \left[ \sup\limits_{t \leq T} \beta^{2}(t) \right] \ \leq \ 2 e^{2T} \left( E \left[ \beta^{2}(0) \right] + 4 \left( g_{T}^{\ast}(0) \right)^{2} E \left[ \Lambda_{T}^{B} \right] + 4 G_{2}^{\ast}(T) \, E \left[ \widehat{\Lambda}_{T}^{H} \right] \right) \ < \ \infty .
	\end{equation}
	
	Recall that the contact rate $\widetilde{\beta} = \widetilde{\beta}_{max} f(\beta)$ is uniquely determined by process $\beta$ -- it is $\mathbb{F}$-adapted and c\`{a}dl\`{a}g, with values in $\left[ \widetilde{\beta}_{max} f_{min}, \widetilde{\beta}_{max} \right]$. Its dynamics follows from the It\^{o} formula for semimartingales (\cite{Protter}, Theorem $32$, page $78$) applied to the function $\widetilde{\beta}_{max} f$, by using $[\beta, \beta]^{c}_{t} = \int_{0}^{t} \sigma_{s}^{2}(0, \beta(s-)) \lambda_{s}^{B} \, ds$ and $\Delta \beta(s) = \sigma_{s}(z, \beta(s-))$ for a jump of size $z$ at time $s$. 
	
	To write the SDE for $\widetilde{\beta}$ in the form \eqref{dfbeta}, the two integrals in it must be well defined. By $\it{(f2), (f3)}$, the second-order Taylor estimate and \eqref{sigma-envelope}, for all $t$ and $z$,
	\begin{equation*} \label{a5-b5-bounds}
		\left| a_{5}(t) \right| \leq \widetilde{\beta}_{max} \left( K_{f^{\prime}} \left| \beta(t) \right|
		+ \frac{1}{2} K_{f^{\prime \prime}} \left( g_{t}^{2}(0) \lambda_{t}^{B} + \lambda_{t}^{H} \int\limits_{\mathbb{R}_{0}} g_{t}^{2}(z) \, \nu(dz) \right) \right), \quad \left| b_{5}(t, z) \right| \leq \widetilde{\beta}_{max} K_{f^{\prime}} \, g_{t}(z),
	\end{equation*}
	and therefore, by \eqref{beta-bound}, \eqref{sigma-envelope-moments} and $\it{(l3)}$, for every $T > 0$
	\begin{align*}
		\int\limits_{0}^{T} \left| a_{5}(t) \right| dt
		& \leq \widetilde{\beta}_{max} \left( K_{f^{\prime}} C_{T}^{(\beta)} T
		+ \frac{1}{2} K_{f^{\prime \prime}} \left( \left( g_{T}^{\ast}(0) \right)^{2} \Lambda_{T}^{B} + G_{2}^{\ast}(T) \, \widehat{\Lambda}_{T}^{H} \right) \right) < \infty, \quad P\text{-a.s.}, \\ 
		E \left[ \int\limits_{0}^{T} \int\limits_{\mathbb{R}} b_{5}^{2}(t, z) \, \Lambda(dt, dz) \right]
		& \leq \widetilde{\beta}_{max}^{2} K_{f^{\prime}}^{2} \left( \left( g_{T}^{\ast}(0) \right)^{2} E \left[ \Lambda_{T}^{B} \right] + G_{2}^{\ast}(T) \, E \left[ \widehat{\Lambda}_{T}^{H} \right] \right) < \infty . \label{b5-square-integrable}
	\end{align*}
	The $dt$-integral in the fifth equation of \eqref{sde1} is thus $P$-a.s. absolutely convergent on bounded time-intervals, and the integral with respect to the random measure $\mu$ is a square-integrable martingale starting from zero. To conclude, the existence and the uniqueness of the contact rate process $\widetilde{\beta}$ are inherited from process $\beta$ and the corresponding SDE \eqref{dfbeta} is well-posed.
	
	\bigskip
	
	\noindent\textbf{The random SIRV model} \\
	It remains to prove the existence, uniqueness and positivity of solution of the first four equations of system \eqref{sde-vector-form}. Let $\Omega_{0} \subset \Omega$, with $P(\Omega_{0}) = 1$, be the set on which the path $t \mapsto \widetilde{\beta}(t, \omega)$ is c\`{a}dl\`{a}g with values in $\left[ \widetilde{\beta}_{max} f_{min}, \widetilde{\beta}_{max} \right]$, and fix $\omega \in \Omega_{0}$. Write $x = \left( x_{1}, x_{2}, x_{3}, x_{4} \right)$ for $\left( S, I, R, V \right)$, $n(x) := x_{1} + x_{2} + x_{3} + x_{4}$ and $D := \langle 0, \infty \rangle^{4}$, and let $F(t, x) := \left( a_{1}, a_{2}, a_{3}, a_{4} \right)(x)$ denote the first four components of the drift \eqref{sde1-vector-form-a} evaluated along the fixed path of $\widetilde{\beta}(t)$ and with $n(x)$ in place of $N(t)$. The compartment equations from system \eqref{sde-vector-form} form the following initial value problem:
	\begin{equation} \label{RODE}
		\frac{d x(t)}{dt} = F(t, x(t)), \qquad x(0) = \left( S(0), I(0), R(0), V(0) \right) \in D.
	\end{equation}
	Equation \eqref{RODE} is the ODE whose right-hand side depends on $t$ only through the fixed path of the contact rate $\widetilde{\beta}$. More precisely, it is a random ODE (RODE) in the sense of \cite{HanKloeden2017} (see also \cite{Arnold1998}, Section $2.2$). 
	
	The function $F$ satisfies the Carath\'{e}odory conditions on $[0, \infty \rangle \times D$, i.e.
	\begin{enumerate}
	\item for each $x \in D$ the map $t \mapsto F(t, x)$ is Borel measurable since it is an affine function of the c\`{a}dl\`{a}g path of $\widetilde{\beta}$,
	\item for each $t \geq 0$ the map $x \mapsto F(t, x)$ is continuously differentiable on $D$,
	\item on every compact set $C \subset D$ both $\left| F(t, x) \right|$ and the Lipschitz constant of $x \mapsto F(t, x)$ on $C$ are bounded uniformly in $t$, since the only non-affine terms are of the form $\widetilde{\beta}(t) \, x_{i} x_{j} / n(x)$, where $\widetilde{\beta}(t) \leq \widetilde{\beta}_{max}$ and $n(x)$ is bounded away from zero on $C$. 
	\end{enumerate}
	Under these conditions the Carath\'{e}odory existence theorem and the uniqueness theorem for locally Lipschitz functions (see e.g. \cite{CoddingtonLevinson1955}, Chapter $2$, or \cite{Hale1980}, Chapter I, Section $5$) yield a unique absolutely continuous solution $x$ of equation \eqref{RODE} on a maximal interval $[0, \tau_{0}(\omega) \rangle$, where $\tau_{0}$ denotes the explosion time. 
	
	Since the initial value $x(0)$ is $\mathcal{F}_{0}$-measurable and $F$ depends on $\omega$ only through the $\mathbb{F}$-adapted path of $\widetilde{\beta}$, the solution at time $t$ is a measurable functional of $\left( x(0), \, \widetilde{\beta}(s), \, 0 \leq s \leq t \right)$. 
	Therefore, $$\left( \left( S(t), I(t), R(t), V(t) \right), \, t \geq 0 \right)$$ is an $\mathbb{F}$-adapted process on $\Omega_{0}$. For the measurability of solutions of RODEs we refer to \cite{HanKloeden2017}, Chapter $2$. Together with already proved existence and uniqueness of solution of the SDE for the contact rate $\widetilde{\beta}$, it follows that the system \eqref{sde-vector-form} admits a unique solution up to the explosion time $\tau_{0}$.
	
	\begin{remark}\label{rem4.1}
	 If the equation \eqref{sde-vector-form} is observed as a SDE driven by the martingale random measure $\mu$, the same conclusion can be obtained without fixing the path of the contact rate $\widetilde{\beta}$. Namely, \eqref{sde-vector-form} is a particular form of the multidimensional SDE observed in \cite{jac2}, where Jacod proved that a solution of a system of this form exists and is unique (Theorem 1), that it is maximal (Theorem 2) and that it does not explode (Theorem 3). Alternatively, the conditions from Cohen \& Elliott \cite{Cohen} (Theorem 16.3.11, page $415$), under which an SDE driven by a semimartingale random measure with uniformly Lipschitz coefficients admits a unique solution, can also be verified directly. Since $b$ does not depend on the state vector and its first four components vanish, only $a_{1}, \ldots, a_{4}$ have to be examined, and all of their terms are linear except those of the form $\widetilde{\beta}(t) I(t) S(t) / N(t)$. However, on the compact domain $\left[ e^{-k}, e^{k} \right]^{4}$, $k \in \mathbb{N}$, whose exit time is the stopping time $\tau_{k}$, the latter are Lipschitz with a constant of the order $\widetilde{\beta}_{max} e^{2k}$. Therefore, this approach  also produces the solution only up to the exit time from that domain, so the localisation that follows is needed in either case. 
	 
	 Also, this is the usual approach if noise is introduced in the first four equations of \eqref{sde-vector-form}, discussed in Remark \ref{rem5.1}. In that case, the components $b_{1}, \ldots, b_{4}$ are then no longer zero, the compartments become semimartingales, and the pathwise construction above is no longer available.
	\end{remark}
	
	\bigskip
	
	\noindent\textbf{Global solution} \\
	Furthermore, we show that the solution of \eqref{sde-vector-form} stays in $\langle 0, \infty \rangle^{4} \times \left[ \widetilde{\beta}_{max} f_{min}, \widetilde{\beta}_{max} \right]$. For every $k \in \mathbb{N}$ we define the stopping time
	\begin{equation*} \label{stopping-time-tau-k}
		\tau_{k} = \inf \left\{ t \in [0, \tau_{0} \rangle: \min{\{ \ln{S(t)}, \ln{I(t)}, \ln{R(t)}, \ln{V(t)} \}} \leq -k \textrm { or } \right.
	\end{equation*}
	\begin{equation*}
		\hspace*{27mm} \left. \max{\{\ln{S(t)}, \ln{I(t)}, \ln{R(t)}, \ln{V(t)} \}} \geq k \right\},
	\end{equation*}
	where the usual agreement $\inf{\emptyset} = \infty$ holds. Note that $\tau_{k}$ increases with increasing $k$ and denote $\tau_{\infty} = \lim\limits_{k \to \infty} \tau_{k}$, $P$-a.s. Since $\tau_{\infty} \leq \tau_{0}$ $P$-a.s., to complete the proof we need to show that $\tau_{\infty} = \infty$ $P$-a.s., meaning that $\left( S(t), I(t), R(t), V(t) \right)$ $P$-a.s. remains in $\langle 0, \infty \rangle^{4}$ for all $t > 0$.
	
	Let us suppose the opposite, i.e. let us suppose that there exists a triple of constants $(k_{0}, T, \varepsilon) \in \mathbb{N} \times \mathbb{R}_{+} \times \langle 0, 1 \rangle$ such that
	\begin{equation} \label{assumption-tauk-T}
		P(\tau_{k} \leq T) \geq \varepsilon, \quad \forall k \geq k_{0}.
	\end{equation}
	First, let us observe the equation \eqref{N-equation} for $N(t)$, obtained by summing the first four equations of system \eqref{sde1}. Then, from the explicit solutions of the differential equations for $I$ and $R$, observed here for $t \leq \tau_{k}$
	\begin{align}
		I(t) & = I(0) \exp{\left\{ \int\limits_{0}^{t} \left( \frac{\widetilde{\beta}(s)}{N(s)} \left( S(s) + \delta V(s) \right) - \left( \kappa_{1} + \theta \right) \right) ds \right\}}, \quad I(0) > 0, \label{I-explicit} \\
		R(t) & = e^{-(\kappa + \gamma) t} \left( R(0) + \theta \int\limits_{0}^{t} I(s) \, e^{(\kappa + \gamma) s} \, ds \right), \quad R(0) > 0, \label{R-explicit}
	\end{align}
	it follows that $I(t) > 0$ and $R(t) > 0$ for all $t \leq \tau_{k}$. Since all four compartments are positive on $[0, \tau_{k}]$, we have that $S(t), I(t), R(t), V(t) < N(t)$ for all $t \leq \tau_{k}$. Therefore, from equation \eqref{N-equation} it follows that
	\begin{equation*}
		\kappa K - \kappa_{1} N(t) \leq \frac{d N(t)}{dt} \leq \kappa (K -  N(t)).
	\end{equation*}
	Gronwall's inequality applied to both sides of the above inequality, or directly solving the equation \eqref{N-equation}, gives \eqref{N-bounds} for all $t \leq \tau_{k}$. In particular, since $0 < N(0) \leq K$ by assumption, \eqref{N-bounds} gives $0 < N(t) \leq K$ for $t \leq \tau_{k}$, and together with $S(t), I(t), R(t), V(t) < N(t)$ it follows that
	\begin{equation} \label{N-bounded-by-K}
		\max \left\{ S(t), I(t), R(t), V(t) \right\} \leq N(t) \leq K, \qquad t \leq \tau_{k} .
	\end{equation}
	Further, for any $t \leq T$, let us define a twice continuously differentiable function $Y \colon \langle 0, \infty \rangle^{4} \times \mathbb{R} \to \mathbb{R}_{+}$:
	\begin{equation*} \label{Lyapunov-function}
		Y(X(t)) = Y \left(X_{1}(t), X_{2}(t), X_{3}(t), X_{4}(t), X_{5}(t) \right) := \sum\limits_{i = 1}^{4} \left( X_{i}(t) - 1 - \ln{X_{i}(t)} \right) + \frac{1}{2} X_{5}^{2}(t),
	\end{equation*}
	where $X_{1}(t) = S(t)$, $X_{2}(t) = I(t)$, $X_{3}(t) = R(t)$, $X_{4}(t) = V(t)$ and $X_{5}(t) = \widetilde{\beta}(t)$. The non-negativity of the function $Y$ follows from the inequality $(x - 1 - \ln{x}) \geq 0$ for all $x > 0$. Note also that $E \left[ Y(X(0)) \right] < \infty$. 
	
	By applying the It\^{o} formula for semimartingales (\cite{Protter}, Theorem $33$, page $81$) to $Y(X(t))$, and by denoting $\psi := \frac{1}{2} \widetilde{\beta}_{max}^{2} f^{2}$, from $(f1)-(f3)$ we have
	\begin{equation} \label{psi-derivatives}
		\left| \psi^{\prime} \right| = \widetilde{\beta}_{max}^{2} \left| f f^{\prime} \right| \leq \widetilde{\beta}_{max}^{2} K_{f^{\prime}} =: C_{1},
		\qquad
		\left| \psi^{\prime \prime} \right| = \widetilde{\beta}_{max}^{2} \left| \left( f^{\prime} \right)^{2} + f f^{\prime \prime} \right| \leq \widetilde{\beta}_{max}^{2} \left( K_{f^{\prime}}^{2} + K_{f^{\prime \prime}} \right) =: 2 C_{2}.
	\end{equation}
	Now it follows that for every $t \leq \tau_{k}$
	\begin{equation} \label{Y-Ito}
		d Y(X(t)) = L(X(t)) \, dt + d \, l(X(t)),
	\end{equation}
	where
	\begin{align}
		L(X(t)) & = \sum\limits_{i=1}^{4} \left( 1 - \frac{1}{X_{i}(t)} \right) \frac{d X_{i}(t)}{dt}
		+ \psi^{\prime}\left( \beta(t-) \right) \beta(t)
		+ \frac{1}{2} \psi^{\prime \prime}\left( \beta(t-) \right) \sigma_{t}^{2}\left( 0, \beta(t-) \right) \lambda_{t}^{B} \nonumber \\
		& \qquad + \lambda_{t}^{H} \int\limits_{\mathbb{R}_{0}} \left( \psi\left( \beta(t-) + \sigma_{t}(z, \beta(t-)) \right) - \psi\left( \beta(t-) \right) - \psi^{\prime}\left( \beta(t-) \right) \sigma_{t}(z, \beta(t-)) \right) \nu(dz) , \label{L-definition} \\
		d \, l(X(t)) & = \int\limits_{\mathbb{R}} \Psi_{t}(z) \, \mu(dt, dz), \label{l-definition}
	\end{align}
	and where
	\begin{equation*}
		\Psi_{t}(z) := \mathbf{1}_{\{ z = 0 \}} \, \psi^{\prime}\left( \beta(t-) \right) \sigma_{t}\left( 0, \beta(t-) \right) + \mathbf{1}_{\{ z \neq 0 \}} \left( \psi\left( \beta(t-) + \sigma_{t}(z, \beta(t-)) \right) - \psi\left( \beta(t-) \right) \right) .
	\end{equation*}
	
	To estimate the sum from $L(X(t))$ in \eqref{L-definition}, recall that on $[0, \tau_{k}]$ all four coordinate processes are positive, hence $N(t) > 0$, $I(t) < N(t)$ and therefore it follows that
	\begin{equation*} \label{transmission-bounds}
		0 \leq \frac{\widetilde{\beta}(t) I(t)}{N(t)} \leq \widetilde{\beta}_{max},
		\qquad
		0 \leq \frac{\widetilde{\beta}(t)}{N(t)} \left( S(t) + \delta V(t) \right) \leq \widetilde{\beta}_{max},
	\end{equation*}
	Discarding the non-positive contributions, and using estimates given above, it follows:
	\begin{align*}
		& \frac{dN}{dt} = \kappa \left( K - N \right) - \left( \kappa_{1} - \kappa \right) I \leq \kappa K, \\
		& - \frac{dS(t)}{S(t)dt} = - \frac{\lambda H}{S} + \left( \rho + \kappa \right) + \frac{\widetilde{\beta} I}{N} - \frac{\alpha V}{S} - \frac{\gamma R}{S} \leq \rho + \kappa + \widetilde{\beta}_{max}, \\
		& - \frac{dI(t)}{I(t)dt} = - \frac{\widetilde{\beta}}{N} \left( S + \delta V \right) + \kappa_{1} + \theta \leq \kappa_{1} + \theta, \\
		& - \frac{dR(t)}{R(t)dt} = - \frac{\theta I}{R} + \kappa + \gamma \leq \kappa + \gamma, \\
		& - \frac{dV(t)}{V(t)dt} = - \frac{\rho S}{V} + \left( \kappa + \alpha \right) + \frac{\delta \widetilde{\beta} I}{N} \leq \kappa + \alpha + \delta \widetilde{\beta}_{max}.
	\end{align*}
	By adding these five estimates, we obtain the following:
	\begin{equation*} \label{LV-compartments}
		\sum\limits_{i=1}^{4} \left( 1 - \frac{1}{X_{i}(t)} \right) \frac{d X_{i}(t)}{dt}
		\ \leq \ \kappa K + \left( 1 + \delta \right) \widetilde{\beta}_{max} + \rho + 3 \kappa + \kappa_{1} + \theta + \gamma + \alpha \ =: \ C_{0} ,
	\end{equation*}
   where $C_{0}$ is independent of both the state variables and of $k$. The remaining three terms of \eqref{L-definition} are estimated by \eqref{psi-derivatives}, by the second-order Taylor estimate $\left| \psi(b + u) - \psi(b) - \psi^{\prime}(b) u \right| \leq C_{2} u^{2}$ and by $(\sigma2)$, giving for $t \leq T$
	\begin{equation*} \label{L-bound}
		L(X(t)) \ \leq \ C_{0} + C_{1} \left| \beta(t) \right| + C_{2} \left[ \left( g_{T}^{\ast}(0) \right)^{2} \lambda_{t}^{B} + G_{2}^{\ast}(T) \lambda_{t}^{H} \right] .
	\end{equation*}
	Consequently, by the Cauchy-Schwarz inequality, relation \eqref{beta-moment} and $(l3)$, it follows:
	\begin{equation} \label{L-expectation}
		E \left[ \int\limits_{0}^{T \wedge \tau_{k}} L(X(t)) \, dt \right]
		\ \leq \ C_{0} T + C_{1} T \sqrt{ E \left[ \sup\limits_{t \leq T} \beta^{2}(t) \right] }
		+ C_{2} \left( \left( g_{T}^{\ast}(0) \right)^{2} E \left[ \Lambda_{T}^{B} \right] + G_{2}^{\ast}(T) \, E \left[ \widehat{\Lambda}_{T}^{H} \right] \right)
		 < \infty.
	\end{equation}
	Note that this bound does not depend on $k$. Further, to estimate \eqref{l-definition} related to $l(X(t))$, as
	$\left| \Psi_{t}(z) \right| \leq C_{1} g_{t}(z)$ for every $z$, by \eqref{psi-derivatives}, the mean value theorem, \eqref{sigma-envelope} and $(l3)$, it follows:
	\begin{equation*} \label{l-square-integrable}
		E \left[ \left( \int\limits_{0}^{T} \int\limits_{\mathbb{R}} \Psi_{t}(z) \, \mu(dt, dz) \right)^{2} \right] = E \left[ \int\limits_{0}^{T} \int\limits_{\mathbb{R}} \Psi_{t}^{2}(z) \, \Lambda(dt, dz) \right]
		\leq C_{1}^{2} \left( \left( g_{T}^{\ast}(0) \right)^{2} E \left[ \Lambda_{T}^{B} \right] + G_{2}^{\ast}(T) \, E \left[ \widehat{\Lambda}_{T}^{H} \right] \right) < \infty.
	\end{equation*}
	Since the random measure $\mu$ is a zero-mean martingale, and the integrand is a square-integrable predictable process, the process $l$ is a square-integrable martingale starting from $0$, and therefore
	\begin{equation*}
		E \left[ \int\limits_{0}^{T \wedge \tau_{k}} \int\limits_{\mathbb{R}} \Psi_{t}(z) \, \mu(dt, dz) \right] = 0 .
	\end{equation*}
	Taking expectations in \eqref{Y-Ito} over $\left[ 0, T \wedge \tau_{k} \right]$ and using \eqref{L-expectation} together with $E \left[ Y(X(0)) \right] < \infty$, we obtain
	\begin{equation} \label{Y-expectation-finite}
		E \left[ Y \left( X(T \wedge \tau_{k}) \right) \right] < \infty , \qquad \text{uniformly in } k.
	\end{equation}
	
	Now let $A_{k} = \left\{ \tau_{k} \leq T \right\}$ for $k \geq k_{0}$. According to the assumption \eqref{assumption-tauk-T}, it follows that $P \left( A_{k} \right) \geq \varepsilon$. Moreover, for every $\omega \in A_{k}$ at least one of the random variables $X_{1}(\tau_{k}) = S(\tau_{k})$, $X_{2}(\tau_{k}) = I(\tau_{k})$, $X_{3}(\tau_{k}) = R(\tau_{k})$, $X_{4}(\tau_{k}) = V(\tau_{k})$ is either less or equal than $e^{-k}$ or greater or equal than $e^{k}$, and since all the summands of $Y$ are non-negative it follows that
	\begin{equation*}
		Y \left( S(\tau_{k} \wedge T), I(\tau_{k} \wedge T), R(\tau_{k} \wedge T), V(\tau_{k} \wedge T), \widetilde{\beta}(\tau_{k} \wedge T) \right)
		\ \geq \ \min{\left\{ e^{k} - 1 - k, \ e^{-k} - 1 + k \right\}} ,
	\end{equation*}
	and therefore
	\begin{equation*} \label{E-Y-tauk-max-T}
		E \left[ Y \left( S(\tau_{k} \wedge T), I(\tau_{k} \wedge T), R(\tau_{k} \wedge T), V(\tau_{k} \wedge T), \widetilde{\beta}(\tau_{k} \wedge T) \right) \right]
		\ \geq \ \varepsilon \, \min{\left\{ e^{k} - 1 - k, \ e^{-k} - 1 + k \right\}} .
	\end{equation*}
	Therefore, if $k \to \infty$, it follows that
	\begin{equation*} \label{thm1-proof-contradiction}
		E \left[ Y \left( S(\tau_{k} \wedge T), I(\tau_{k} \wedge T), R(\tau_{k} \wedge T), V(\tau_{k} \wedge T), \widetilde{\beta}(\tau_{k} \wedge T) \right) \right] = \infty,
	\end{equation*}
	which is a contradiction with \eqref{Y-expectation-finite}, meaning that the assumption \eqref{assumption-tauk-T} is false. In particular, it follows that $\tau_{\infty} = \infty$ $P$-a.s., and therefore $\tau_{0} = \infty$. It follows that the solution is global and remains in $\langle 0, \infty \rangle^{4} \times \left[ \widetilde{\beta}_{max} f_{min}, \widetilde{\beta}_{max} \right]$ for all $t \geq 0$.
	
	Furthermore, since \eqref{N-bounds} and \eqref{N-bounded-by-K} were established for $t \leq \tau_{k}$ and $\tau_{k} \to \infty$, they hold for all $t \geq 0$. By letting $t \to \infty$ in \eqref{N-bounds}, we obtain \eqref{N-boundaries}.
	
	\bigskip
	
	\noindent\textbf{The positively invariant set} \\
	Let $X(0) \in \Gamma$ $P$-a.s., i.e. let \eqref{initial-condition} hold with $\frac{\kappa}{\kappa_{1}} K \leq N(0) \leq K$. By claims (i) and (ii) of Theorem \ref{thm-global-positive-solution}, it follows that $X(t) \in \langle 0, K]^{4} \times \left[ \widetilde{\beta}_{max} f_{min}, \widetilde{\beta}_{max} \right]$ for all $t \geq 0$, and by \eqref{N-bounds} the population size satisfies $\frac{\kappa}{\kappa_{1}} K \leq N(t) \leq K$ for all $t \geq 0$. Hence, $X(t) \in \Gamma$ for all $t \geq 0$ $P$-a.s., i.e. the set $\Gamma$ given in \eqref{gamma-set} is the positively invariant set of the process $X$.
\end{proof}

\begin{remark} \label{rem-finite-activity}
	If $\nu(\mathbb{R}_{0}) < \infty$, which is the finite-activity case not assumed in this paper, the proof of Theorem \ref{thm-global-positive-solution} simplifies at one point. Therefore, the process $\beta$ has finitely many jumps on every bounded time-interval and the contact rate $\widetilde{\beta}$ is continuous between consecutive jump times. Then, the construction of the solution of equation \eqref{RODE} can be carried out with the classical Picard-Lindel\"{o}f theorem on each inter-jump interval, the solution being restarted at every jump time from its current value. The rest of the proof stays the same, as well as the proofs of Theorems \ref{thm-extinction} and \ref{thm-persistence}.
\end{remark}

	\section{Extinction} \label{sec:extinction}

	In this section we generate the conditions under which the disease dies out from the population (but the population survives). In this analysis $t \to \infty$, an therefore it requires three strengthenings of assumptions of Section \ref{sec:EUPsolution}, which were all stated on bounded time-horizons. The additional requirements are the following, and they are used only in this section and in Section \ref{sec:persistence}:
	\begin{itemize}
		\item[$\it{(f4)}$] the contact rate saturates, in the sense that
		\begin{equation*} \label{f4-assumption}
			\sup\limits_{x \in \mathbb{R}} \left| x f^{\prime}(x) \right| =: K_{x f^{\prime}} < \infty,
		\end{equation*}
		\item[$(\it{\sigma2^{\ast})}$] the function $g$ assumed to exist in $\it{(\sigma2)}$ is bounded and has finite second and fourth moments uniformly in time, i.e.
		\begin{equation} \label{sigma-envelope-global}
			g^{\ast}(0) := \sup\limits_{t \geq 0} g_{t}(0) < \infty,
			\quad
			G_{2} := \sup\limits_{t \geq 0} \int\limits_{\mathbb{R}_{0}} g_{t}^{2}(z) \, \nu(dz) < \infty,
			\quad
			G_{4} := \sup\limits_{t \geq 0} \int\limits_{\mathbb{R}_{0}} g_{t}^{4}(z) \, \nu(dz) < \infty  \quad P\text{-a.s.},
		\end{equation}
		\item[$\it{(l4)}$] the stochastic clock has a finite long-run average speed, i.e. for $l = B, H$,
		\begin{equation*} \label{l4-assumption}
			\limsup\limits_{t \to \infty} \frac{1}{t} \int\limits_{0}^{t} \lambda^{l}_{s} \, ds =: \ell^{l} < \infty, \quad P\text{-a.s.}
		\end{equation*}
	\end{itemize} 
    Note that $\it{(\sigma2^{\ast})}$ implies $\it{(\sigma2)}$, since $g_{T}^{\ast}(0) \leq g^{\ast}(0)$ and $G_{2}^{\ast}(T) \leq G_{2}$ for every $T$, so replacing $\it{(\sigma2)}$ by $\it{(\sigma2^{\ast})}$ changes nothing in Section \ref{sec:EUPsolution}. In $\it{(l4)}$, $\frac{1}{t} \int_{0}^{t} \lambda_{s}^{l} \, ds$ is the average speed of the stochastic clock over $[0, t]$, and $\ell^{l}$ is the long-run rate at which the operational time accumulates with respect to the calendar time. Condition $\it{(l4)}$ requires this rate to be finite, but not to vanish (e.g. in case of constant intensities, $\lambda^{B} = c_{B}$ and $\lambda^{H} = c_{H}$, we have $\ell^{B} = c_{B}$ and $\ell^{H} = c_{H}$).

	Two long-run averages of the jump coefficient appear in the threshold for the extinction of the disease:
	\begin{equation} \label{ell-definitions}
		\ell_{B} := \limsup\limits_{t \to \infty} \frac{1}{t} \int\limits_{0}^{t} \sigma_{s}^{2} \left( 0, \beta(s-) \right) \lambda_{s}^{B} \, ds,
		\quad
		\ell_{H} := \limsup\limits_{t \to \infty} \frac{1}{t} \int\limits_{0}^{t} \int\limits_{\mathbb{R}_{0}} \sigma_{s}^{2} \left( z, \beta(s-) \right) \lambda_{s}^{H} \, \nu(dz) \, ds, \quad P\text{-a.s.},
	\end{equation}
	and note the distinction from the averages $\ell^{B}, \ell^{H}$ of the stochastic clock. In epidemiological terms, $\ell_{B}$ and $\ell_{H}$ measure, on average per unit of calendar time, how strongly the underlying process $\beta$ is moved by the continuous behavioural fluctuations and by the shocks, respectively. By \eqref{sigma-envelope} and \eqref{sigma-envelope-global}, $\sigma_{s}^{2}(0, \cdot) \leq \left( g^{\ast}(0) \right)^{2}$ and $\int_{\mathbb{R}_{0}} \sigma_{s}^{2}(z, \cdot) \, \nu(dz) \leq G_{2}$ for all $s \geq 0$, so by $(l4)$
	\begin{equation*} \label{ell-bounds}
		\ell_{B} \leq \left( g^{\ast}(0) \right)^{2} \ell^{B} < \infty,
		\qquad
		\ell_{H} \leq G_{2} \, \ell^{H} < \infty, \quad P\text{-a.s.},
	\end{equation*}
	i.e. $\ell_{B}$ and $\ell_{H}$ are $P$-a.s. finite random variables.

	\begin{theorem} \label{thm-extinction}
		Let $\it{(f4)}$, $\it{(\sigma2^{\ast})}$ and $\it{(l4)}$ hold, and let the initial value be $\mathcal{F}_{0}$-measurable with values in the positively invariant set $\Gamma$ given in \eqref{gamma-set}. If the long-run averages $\ell_{B}$ and $\ell_{H}$ given by \eqref{ell-definitions} satisfy
		\begin{equation} \label{extinction-condition}
			\ell_{B} + \ell_{H} \ < \ \frac{2 \left( \kappa_{1} + \theta - \widetilde{\beta}_{max} - K_{x f^{\prime}} \right)}{K_{f^{\prime}}^{2} + K_{f^{\prime \prime}}}, \quad P\text{-a.s.},
		\end{equation}
		then
		\begin{equation*}
			\lim\limits_{t \to \infty} I(t) = \lim\limits_{t \to \infty} R(t) = 0 \quad \text{and} \quad \limsup\limits_{t \to \infty} \left( S(t) + V(t) \right) \leq K, \qquad P\text{-a.s.}
		\end{equation*}
	\end{theorem}

	\begin{proof}
		The proof rests on a strong law of large numbers for local martingales, where assumption $\it{(l4)}$ is needed. The time-change processes $\Lambda^{B}_{t}$ and $\widehat{\Lambda}^{H}_{t}$ are non-decreasing and $P$-a.s. finite by $\it{(l1)}$ and $\it{(l3)}$, so by $\it{(l4)}$ there is a $P$-a.s. finite random variable $C$ such that $\Lambda^{B}_{t} + \widehat{\Lambda}^{H}_{t} \leq C (1 + t)$ for all $t \geq 0$. In general, let $M$ be a locally square-integrable martingale with $M(0) = 0$, whose predictable quadratic variation satisfies $\langle M, M \rangle_{t} \leq c \left( \Lambda^{B}_{t} + \widehat{\Lambda}^{H}_{t} \right)$ for all $t \geq 0$ and some $P$-a.s. finite random variable $c$. Integrating by parts, for every $T > 0$
		\begin{equation*}
			\int\limits_{0}^{T} \frac{d \langle M, M \rangle_{s}}{(1+s)^{2}}
			= \frac{\langle M, M \rangle_{T}}{(1+T)^{2}} + 2 \int\limits_{0}^{T} \frac{\langle M, M \rangle_{s}}{(1+s)^{3}} \, ds
			\leq \frac{c C}{1+T} + 2 c C \int\limits_{0}^{T} \frac{ds}{(1+s)^{2}} \leq 3 c C ,
		\end{equation*}
		so that $\int_{0}^{\infty} (1+s)^{-2} \, d \langle M, M \rangle_{s} < \infty$ $P$-a.s. Then the strong law of large numbers for local martingales (\cite{Lipster}, Chapter 2, Theorem $10$, page $144$; see also \cite{Lipster-paper}, Theorem $1$) gives
		\begin{equation} \label{SLLN-consequence}
			\lim\limits_{t \to \infty} \frac{M(t)}{t} = 0, \quad P\text{-a.s.}
		\end{equation}
		
		We now turn to the conditions under which the disease dies out. Consider the process $F(I(t), \widetilde{\beta}(t)) = \ln{I(t)} + \frac{1}{2} f^{2} \left( \beta(t) \right)$, where
		\begin{equation*} \label{extinction-function-F}
			F(x, y) := \ln{(x)} + \frac{y^{2}}{2 \widetilde{\beta}_{max}^{2}}, \quad x > 0, \quad y \in \mathbb{R}.
		\end{equation*}
		By Theorem \ref{thm-global-positive-solution}, $I$ is positive and absolutely continuous and $\frac{d}{dt} \ln{I(t)} = \frac{\widetilde{\beta}(t)}{N(t)} \left( S(t) + \delta V(t) \right) - \left( \kappa_{1} + \theta \right)$. The It\^{o} formula for semimartingales (\cite{Protter}, Theorem $32$, page $78$), applied to the $C^{2}$ function $\frac{1}{2} f^{2}$ and to $\beta$ given by \eqref{SDE-beta}, with $[\beta, \beta]^{c}_{t} = \int_{0}^{t} \sigma_{s}^{2}(0, \beta(s-)) \lambda_{s}^{B} \, ds$, therefore gives
		\begin{align} \label{extinction-proof-Ito-function-F}
			& F \left( I(t), \widetilde{\beta}(t) \right) = F \left( I(0), \widetilde{\beta}(0) \right) + \int\limits_{0}^{t} \left[ \frac{\widetilde{\beta}(s)}{N(s)} \left( S(s) + \delta V(s) \right) - \left( \kappa_{1} + \theta \right) + f(\beta(s-)) f^{\prime}(\beta(s-)) \beta(s) \right] ds \nonumber \\
			& \qquad + \frac{1}{2} \int\limits_{0}^{t} \left[ \left( f^{\prime}(\beta(s-)) \right)^{2} + f(\beta(s-)) f^{\prime \prime}(\beta(s-)) \right] \sigma_{s}^{2} \left( 0, \beta(s-) \right) \lambda_{s}^{B} \, ds \nonumber \\
			& \qquad + \int\limits_{0}^{t} \int\limits_{\mathbb{R}} f(\beta(s-)) f^{\prime}(\beta(s-)) \sigma_{s} \left( z, \beta(s-) \right) \, \mu(ds, dz) \nonumber \\
			& \qquad + \sum\limits_{0 < s \leq t} \left[ \frac{1}{2} f^{2} \left( \beta(s) \right) - \frac{1}{2} f^{2} \left( \beta(s-) \right) - f \left( \beta(s-) \right) f^{\prime} \left( \beta(s-) \right) \Delta \beta(s) \right].
		\end{align}
		
		Next, we bound the terms of \eqref{extinction-proof-Ito-function-F} one by one. Since $\delta \in [0, 1]$ and $I, R \geq 0$, we have $S + \delta V \leq N$, and therefore by $\it{(f1)}$ it follows:
		\begin{equation*} \label{beta-tilde-sharp-bound}
			\frac{\widetilde{\beta}(t)}{N(t)} \left( S(t) + \delta V(t) \right) \leq \widetilde{\beta}(t) = \widetilde{\beta}_{max} f \left( \beta(t) \right) \leq \widetilde{\beta}_{max}, \quad t \geq 0 .
		\end{equation*}
		Since $\beta(s-) = \beta(s)$ outside the countable set of jump times, $(f1)$ and $(f4)$ give $\left| f(\beta(s-)) f^{\prime}(\beta(s-)) \beta(s) \right| \leq \left| \beta(s) f^{\prime}(\beta(s)) \right| \leq K_{x f^{\prime}}$. With $C_{f} := \frac{1}{2} \left( K_{f^{\prime}}^{2} + K_{f^{\prime \prime}} \right)$, assumptions $(f1)-(f3)$ give $\frac{1}{2} \left| (f^{\prime})^{2} + f f^{\prime \prime} \right| \leq C_{f}$, and Taylor's theorem with the Lagrange form of the remainder, applied to $\frac{1}{2} f^{2}$, bounds each summand of the jump sum by $C_{f} \left( \Delta \beta(s) \right)^{2}$. Discarding the non-negative term $\frac{1}{2} f^{2}(\beta(t))$ on the left-hand side, we obtain
		\begin{align}
			\ln{I(t)} & \leq F \left( I(0), \widetilde{\beta}(0) \right) + \left[ \widetilde{\beta}_{max} - \left( \kappa_{1} + \theta \right) + K_{x f^{\prime}} \right] t \nonumber \\
			& \qquad + \frac{1}{2} \left( K_{f^{\prime}}^{2} + K_{f^{\prime \prime}} \right) \left( \int\limits_{0}^{t} \sigma_{s}^{2} \left( 0, \beta(s-) \right) \lambda_{s}^{B} \, ds + \sum\limits_{0 < s \leq t} \left( \Delta \beta(s) \right)^{2} \right) + M^{f}(t), \label{extinction-master-estimate}
		\end{align}
		where $M^{f}(t) := \int_{0}^{t} \int_{\mathbb{R}} f(\beta(s-)) f^{\prime}(\beta(s-)) \sigma_{s} \left( z, \beta(s-) \right) \mu(ds, dz)$.
		
		Now we divide \eqref{extinction-master-estimate} by $t$ and let $t \to \infty$. Since $(f f^{\prime})^{2} \leq K_{f^{\prime}}^{2}$, the bounds \eqref{sigma-envelope} and \eqref{sigma-envelope-global} give
		\begin{equation*}
			\langle M^{f}, M^{f} \rangle_{t} = \int\limits_{0}^{t} \int\limits_{\mathbb{R}} \left( f(\beta(s-)) f^{\prime}(\beta(s-)) \right)^{2} \sigma_{s}^{2} \left( z, \beta(s-) \right) \Lambda(ds, dz)
			\leq K_{f^{\prime}}^{2} \left[ \left( g^{\ast}(0) \right)^{2} \Lambda^{B}_{t} + G_{2} \, \widehat{\Lambda}^{H}_{t} \right],
		\end{equation*}
		so that $M^{f}(t)/t \to 0$ $P$-a.s. by the argument \eqref{SLLN-consequence}. The jumps of process $\beta$ are of the form  $\Delta \beta(s) = \int_{\mathbb{R}_{0}} \sigma_{s}(z, \beta(s-)) \, H(\{s\}, dz)$, and therefore
		\begin{equation*}
			\sum\limits_{0 < s \leq t} \left( \Delta \beta(s) \right)^{2} = \int\limits_{0}^{t} \int\limits_{\mathbb{R}_{0}} \sigma_{s}^{2} \left( z, \beta(s-) \right) H(ds, dz) = \int\limits_{0}^{t} \int\limits_{\mathbb{R}_{0}} \sigma_{s}^{2} \left( z, \beta(s-) \right) \lambda_{s}^{H} \nu(dz) \, ds + \widetilde{M}(t),
		\end{equation*}
		where $\widetilde{M}(t) := \int_{0}^{t} \int_{\mathbb{R}_{0}} \sigma_{s}^{2}(z, \beta(s-)) \, \widetilde{H}(ds, dz)$ is a locally square-integrable martingale with predictable quadratic variation $\langle \widetilde{M}, \widetilde{M} \rangle_{t} = \int_{0}^{t} \int_{\mathbb{R}_{0}} \sigma_{s}^{4}(z, \beta(s-)) \lambda_{s}^{H} \nu(dz) \, ds \leq G_{4} \, \widehat{\Lambda}^{H}_{t}$, finite by assumption $\it{(\sigma2^{\ast})}$. Hence $\widetilde{M}(t)/t \to 0$ $P$-a.s. by \eqref{SLLN-consequence}, and
		\begin{equation} \label{jump-sum-limit}
			\limsup\limits_{t \to \infty} \frac{1}{t} \sum\limits_{0 < s \leq t} \left( \Delta \beta(s) \right)^{2} = \limsup\limits_{t \to \infty} \frac{1}{t} \int\limits_{0}^{t} \int\limits_{\mathbb{R}_{0}} \sigma_{s}^{2} \left( z, \beta(s-) \right) \lambda_{s}^{H} \nu(dz) \, ds = \ell_{H}, \quad P\text{-a.s.}
		\end{equation}
		Since $F( I(0), \widetilde{\beta}(0))/t \to 0$, it follows that
		\begin{equation} \label{extinction-lnI-limsup}
			\limsup\limits_{t \to \infty} \frac{\ln{I(t)}}{t} \leq \widetilde{\beta}_{max} - \left( \kappa_{1} + \theta \right) + K_{x f^{\prime}} + \frac{1}{2} \left( K_{f^{\prime}}^{2} + K_{f^{\prime \prime}} \right) \left( \ell_{B} + \ell_{H} \right), \quad P\text{-a.s.}
		\end{equation}
		Therefore, if
		\begin{equation*}
			\kappa_{1} + \theta > \widetilde{\beta}_{max} + K_{x f^{\prime}} + \frac{1}{2} \left( K_{f^{\prime}}^{2} + K_{f^{\prime \prime}} \right) \left( \ell_{B} + \ell_{H} \right),
		\end{equation*}
		the right-hand side of inequality \eqref{extinction-lnI-limsup} is strictly negative, $I(t)$ decays exponentially fast as $t \to \infty$ and
		\begin{equation*} \label{extinction-I}
			\lim\limits_{t \to \infty} I(t) = 0, \quad P\text{-a.s.}
		\end{equation*}
		For the recovered class, the solution \eqref{R-explicit} of the equation for $R$ gives:
		\begin{equation*}
			R(t) = R(0) \, e^{-(\kappa + \gamma) t} + \theta \, \frac{\int_{0}^{t} I(s) \, e^{(\kappa + \gamma) s} \, ds}{e^{(\kappa + \gamma) t}}, \quad R(0) > 0.
		\end{equation*}
		The first term vanishes as $t \to \infty$. In the second, the denominator tends to infinity and the numerator is continuously differentiable, and since $\lim\limits_{t \to \infty} I(t) = 0$, $P$-a.s., by L'Hospital's rule it follows that
		\begin{equation*}
			\lim\limits_{t \to \infty} \theta \, \frac{\int_{0}^{t} I(s) \, e^{(\kappa + \gamma) s} \, ds}{e^{(\kappa + \gamma) t}} = \lim\limits_{t \to \infty} \frac{\theta \, I(t) \, e^{(\kappa + \gamma) t}}{(\kappa + \gamma) \, e^{(\kappa + \gamma) t}} = \frac{\theta}{\kappa + \gamma} \lim\limits_{t \to \infty} I(t) = 0, \quad P\text{-a.s.}
		\end{equation*}
		Therefore, $\lim\limits_{t \to \infty} R(t) = 0$, $P$-a.s. Finally, $S(t) + V(t) = N(t) - I(t) - R(t)$, and $\limsup_{t \to \infty} \left( S(t) + V(t) \right) = \limsup_{t \to \infty} N(t) \leq K$ $P$-a.s.
	\end{proof}

\begin{figure}[htbp]
	\centering
	\includegraphics[width=\textwidth]{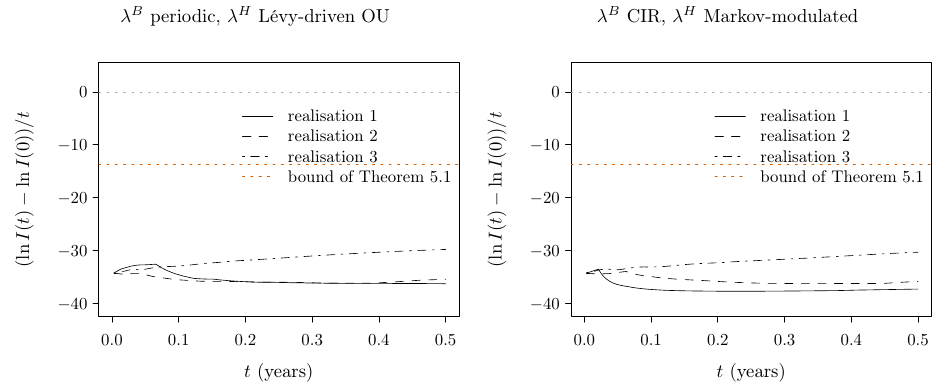}
	\caption{Extinction regime ($\overline{\mathcal{R}}_{0} = 0.7$ from Section~\ref{sec:simulation}): the condition $(\ln I(t) - \ln I(0))/t$ for the three realisations of the noise from Figure~\ref{fig:beta} (solid, dashed and dash-dotted lines) and the two stochastic clocks, together with the bound \eqref{extinction-lnI-limsup} from Theorem~\ref{thm-extinction}.} \label{fig:extinction-rate}
\end{figure}

\begin{figure}[htbp]
	\centering
	\includegraphics[width=\textwidth]{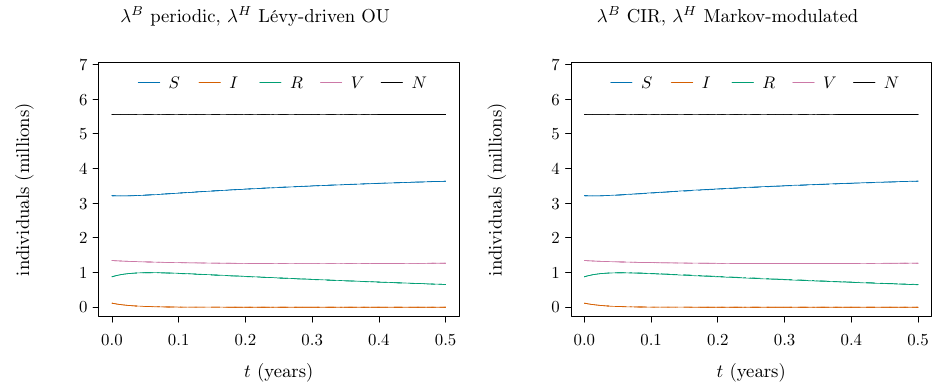}
	\caption{Extinction regime ($\overline{\mathcal{R}}_{0} = 0.7$ from Section~\ref{sec:simulation}): the compartments $S$, $I$, $R$, $V$ and the population size $N$ for the three realisations of the noise from Figure~\ref{fig:beta} (solid, dashed and dash-dotted lines) and the stochastic clocks.} \label{fig:extinction-compartments}
\end{figure}

\begin{figure}[htbp]
	\centering
	\includegraphics[width=\textwidth]{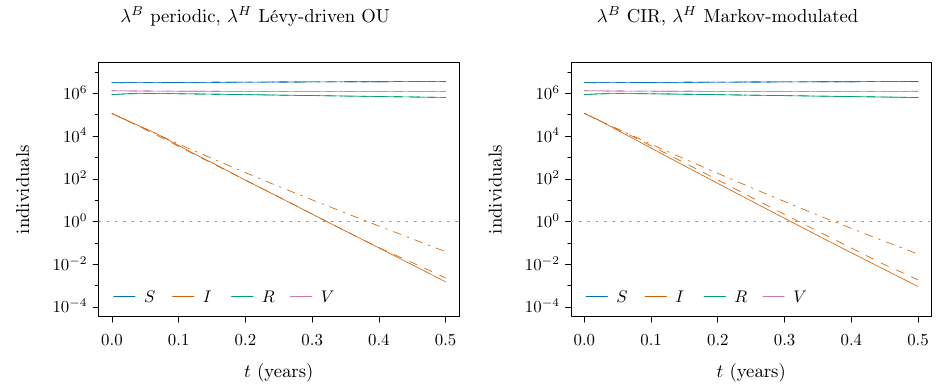}
	\caption{Extinction regime ($\overline{\mathcal{R}}_{0} = 0.7$ from Section~\ref{sec:simulation}): the compartments from Figure~\ref{fig:extinction-compartments} on a logarithmic scale, where the dotted line marks one individual.} \label{fig:extinction-log}
\end{figure}

	\begin{remark} \label{rem-extinction-interpretation}
		Assumption \eqref{extinction-condition} from Theorem \ref{thm-extinction}, under which the disease dies out, can be written in terms of the random basic reproduction number from Subsection \ref{subsec-stochastic-SIRV}, i.e. in terms of its upper bound $\overline{\mathcal{R}}_{0} = \widetilde{\beta}_{max} / (\kappa_{1} + \theta)$:
		\begin{equation*} \label{extinction-condition-R0}
			\overline{\mathcal{R}}_{0} + \frac{1}{\kappa_{1} + \theta} \left( K_{x f^{\prime}} + \frac{K_{f^{\prime}}^{2} + K_{f^{\prime \prime}}}{2} \left( \ell_{B} + \ell_{H} \right) \right) < 1 .
		\end{equation*}
		Since the terms added to $\overline{\mathcal{R}}_{0}$ are non-negative and $\mathcal{R}_{0}(t) \leq \overline{\mathcal{R}}_{0}$, condition \eqref{extinction-condition} implies that $\sup_{t \geq 0} \mathcal{R}_{0}(t) < 1$. In epidemiological terms, the disease dies out if the random basic reproduction number $\mathcal{R}_{0}$ stays below $1$ even at the maximal level of the contact rate, and still leaves room for the noise. For the parameter values from Section~\ref{sec:simulation}, condition \eqref{extinction-condition} holds (see Table~\ref{tab:latent}), and Figure~\ref{fig:extinction-rate} shows that $(\ln I(t) - \ln I(0))/t$ stays below the bound \eqref{extinction-lnI-limsup}, while Figures~\ref{fig:extinction-compartments} and~\ref{fig:extinction-log} show the corresponding paths of the compartments.
	\end{remark}

	\section{Persistence in mean} \label{sec:persistence}
	
	In this section we give conditions under which the disease remains endemic on average, i.e. the time average of the number of infected individuals stays above zero.

	\begin{definition} \label{persistence}
		The random system \eqref{sde1} is said to be persistent in mean if
		\begin{equation*}
			\liminf_{t \to \infty} \frac{1}{t} \int\limits_{0}^{t} I(s) \, ds > 0, \quad P\text{-a.s.},
		\end{equation*}
		for any initial value from the positively invariant set $\Gamma$ given in \eqref{gamma-set}.
	\end{definition}

	The following theorem gives a sufficient condition for persistence in mean, together with an explicit lower bound on the long-run average of infected class $I$.

	\begin{theorem} \label{thm-persistence}
		Let $\it{(f4)}$, $\it{(\sigma2^{\ast})}$ and $\it{(l4)}$ hold, let $\delta > 0$, and let the initial value be $\mathcal{F}_{0}$-measurable with values in the positively invariant set $\Gamma$ given in \eqref{gamma-set}. If the long-run averages $\ell_{B}$ and $\ell_{H}$ given by \eqref{ell-definitions} satisfy
		\begin{equation} \label{persistence-condition}
			\ell_{B} + \ell_{H} \ < \ \frac{2 \left( 1 + f_{min} \right) \left( \left( 1 + f_{min} \right) \left( \delta \widetilde{\beta}_{max} f_{min} - \kappa_{1} - \theta \right) - K_{x f^{\prime}} \right)}{K_{f^{\prime}}^{2} + K_{f^{\prime \prime}} \left( 1 + f_{min} \right)}, \quad P\text{-a.s.},
		\end{equation}
		then the system \eqref{sde1} is persistent in mean. More precisely,
		\begin{equation*} \label{persistence-bound}
			\liminf\limits_{t \to \infty} \frac{1}{t} \int\limits_{0}^{t} I(s) \, ds
			\ \geq \ \frac{\kappa + \gamma}{\kappa + \gamma + \theta} \cdot \frac{\kappa K}{\kappa_{1}}
			\left( 1 - \frac{\Theta}{\delta \widetilde{\beta}_{max} f_{min}} \right) \ > \ 0, \quad P\text{-a.s.},
		\end{equation*}
		where
		\begin{equation} \label{persistence-threshold}
			\Theta := \kappa_{1} + \theta + \frac{K_{x f^{\prime}}}{1 + f_{min}} + \frac{K_{f^{\prime}}^{2} + K_{f^{\prime \prime}} \left( 1 + f_{min} \right)}{2 \left( 1 + f_{min} \right)^{2}} \left( \ell_{B} + \ell_{H} \right).
		\end{equation}
	\end{theorem}
	
			\begin{proof}
			Let us define a function
			\begin{equation*} \label{persistence-function-G}
				G(x, b) := \ln{(x)} - \ln{\left( 1 + f(b) \right)}, \qquad x > 0, \quad b \in \mathbb{R}.
			\end{equation*}
			The first two derivatives of $b \mapsto \ln \left( 1 + f(b) \right)$ are $f^{\prime} / (1 + f)$ and $f^{\prime \prime} / (1 + f) - (f^{\prime})^{2} / (1 + f)^{2}$, and since $f \geq f_{min}$, assumptions $\it{(f2)-(f4)}$ for all $b \in \mathbb{R}$ give
			\begin{align} \label{logf-bounds}
				& \frac{\left| f^{\prime}(b) \right|}{1 + f(b)} \leq \frac{K_{f^{\prime}}}{1 + f_{min}},
				\qquad
				\frac{\left| b f^{\prime}(b) \right|}{1 + f(b)} \leq \frac{K_{x f^{\prime}}}{1 + f_{min}}, \nonumber \\
				& \left| \frac{f^{\prime \prime}(b)}{1 + f(b)} - \frac{\left( f^{\prime}(b) \right)^{2}}{\left( 1 + f(b) \right)^{2}} \right| \leq \frac{K_{f^{\prime \prime}}}{1 + f_{min}} + \frac{K_{f^{\prime}}^{2}}{\left( 1 + f_{min} \right)^{2}} = \frac{K_{f^{\prime}}^{2} + K_{f^{\prime \prime}} \left( 1 + f_{min} \right)}{\left( 1 + f_{min} \right)^{2}} =: 2 C_{G} .
			\end{align}
			By Theorem \ref{thm-global-positive-solution}, $0 < I(t) \leq N(t) \leq K$, so that
			\begin{equation} \label{G-bounded-above}
				G \left( I(t), \beta(t) \right) \leq \ln{K} - \ln{\left( 1 + f_{min} \right)}, \quad t \geq 0, \qquad \text{and hence} \qquad \limsup\limits_{t \to \infty} \frac{G \left( I(t), \beta(t) \right)}{t} \leq 0 .
			\end{equation}
			The It\^{o} formula for semimartingales (\cite{Protter}, Theorem $32$, page $78$) applied to the function $\ln \left( 1 + f(\beta) \right)$, where process $\beta$ is given by \eqref{SDE-beta}, 
			gives the following:
			\begin{align} \label{persistence-Ito}
				G \left( I(t), \beta(t) \right) & = G \left( I(0), \beta(0) \right)
				+ \int\limits_{0}^{t} \left[ \frac{\widetilde{\beta}(s)}{N(s)} \left( S(s) + \delta V(s) \right) - \left( \kappa_{1} + \theta \right) - \frac{f^{\prime}\left( \beta(s-) \right)}{1 + f\left( \beta(s-) \right)} \, \beta(s) \right] ds \nonumber \\
				& \qquad - \frac{1}{2} \int\limits_{0}^{t} \left[ \frac{f^{\prime \prime}\left( \beta(s-) \right)}{1 + f\left( \beta(s-) \right)} - \frac{\left( f^{\prime}\left( \beta(s-) \right) \right)^{2}}{\left( 1 + f\left( \beta(s-) \right) \right)^{2}} \right] \sigma_{s}^{2}\left( 0, \beta(s-) \right) \lambda_{s}^{B} \, ds
				- M^{G}(t) \nonumber \\
				& \qquad - \sum\limits_{0 < s \leq t} \left[ \ln{\frac{1 + f\left( \beta(s) \right)}{1 + f\left( \beta(s-) \right)}} - \frac{f^{\prime}\left( \beta(s-) \right)}{1 + f\left( \beta(s-) \right)} \, \Delta \beta(s) \right],
			\end{align}
			where $M^{G}(t) := \int_{0}^{t} \int_{\mathbb{R}} \frac{f^{\prime}\left( \beta(s-) \right)}{1 + f\left( \beta(s-) \right)} \, \sigma_{s}\left( z, \beta(s-) \right) \mu(ds, dz)$. Next, we bound the terms of expression \eqref{persistence-Ito} from below one by one. 
			
			Since $\beta(s-) = \beta(s)$ outside the countable set of jump times, the bound \eqref{logf-bounds} gives
			\begin{equation} \label{persistence-drift-bound}
				- \frac{f^{\prime}\left( \beta(s-) \right)}{1 + f\left( \beta(s-) \right)} \, \beta(s) \ \geq \ - \frac{\left| \beta(s) f^{\prime}(\beta(s)) \right|}{1 + f(\beta(s))} \ \geq \ - \frac{K_{x f^{\prime}}}{1 + f_{min}} .
			\end{equation}
			The absolute value of the second derivative of $\ln (1 + f)$ is by \eqref{logf-bounds} bounded by the constant $2 C_{G}$, and by Taylor's theorem with the Lagrange form of the remainder, $$\left| \ln \frac{1 + f(b + u)}{1 + f(b)} - \frac{f^{\prime}(b)}{1 + f(b)} \, u \right| \leq C_{G} u^{2}, \quad \forall b, u \in \mathbb{R}.$$ Therefore, the second-order term and the jump term in expression \eqref{persistence-Ito} are bounded from below by
			\begin{equation} \label{persistence-second-order-bounds}
				- C_{G} \int\limits_{0}^{t} \sigma_{s}^{2}\left( 0, \beta(s-) \right) \lambda_{s}^{B} \, ds
				\qquad \text{and} \qquad
				- C_{G} \sum\limits_{0 < s \leq t} \left( \Delta \beta(s) \right)^{2},
			\end{equation}
			respectively. Regarding the incidence term, since $\delta \in \langle 0, 1]$ by assumption of the theorem and $S \geq 0$, we have that $S + \delta V \geq \delta (S + V) = \delta \left( N - I - R \right)$, so that
			\begin{equation*}
				\frac{\widetilde{\beta}(s)}{N(s)} \left( S(s) + \delta V(s) \right) \ \geq \ \delta \, \widetilde{\beta}(s) \left( 1 - \frac{I(s) + R(s)}{N(s)} \right) .
			\end{equation*}
			The second factor in the above upper bound is non-negative, so it may be multiplied by the lower bound $\widetilde{\beta} \geq \widetilde{\beta}_{max} f_{min}$. Since the initial value belongs to the positively invariant set $\Gamma$, $N(s) \geq \kappa K / \kappa_{1}$ for all $s \geq 0$, and therefore
			\begin{equation} \label{incidence-lower}
				\frac{\widetilde{\beta}(s)}{N(s)} \left( S(s) + \delta V(s) \right) \ \geq \ \delta \widetilde{\beta}_{max} f_{min} - \frac{\delta \widetilde{\beta}_{max} \kappa_{1} f_{min}}{\kappa K} \left( I(s) + R(s) \right), \qquad \forall s \geq 0 .
			\end{equation}
			Inserting estimates \eqref{persistence-drift-bound}, \eqref{persistence-second-order-bounds} and \eqref{incidence-lower} in expression \eqref{persistence-Ito} yields the following:
			\begin{align}
				G \left( I(t), \beta(t) \right) & \geq G \left( I(0), \beta(0) \right)
				+ \left( \delta \widetilde{\beta}_{max} f_{min} - \left( \kappa_{1} + \theta \right) - \frac{K_{x f^{\prime}}}{1 + f_{min}} \right) t
				- \frac{\delta \widetilde{\beta}_{max} \kappa_{1} f_{min}}{\kappa K} \int\limits_{0}^{t} \left( I(s) + R(s) \right) ds \nonumber \\
				& \qquad - C_{G} \left( \int\limits_{0}^{t} \sigma_{s}^{2}\left( 0, \beta(s-) \right) \lambda_{s}^{B} \, ds + \sum\limits_{0 < s \leq t} \left( \Delta \beta(s) \right)^{2} \right) - M^{G}(t) . \label{persistence-master-estimate}
			\end{align}
			Since by Theorem \ref{thm-global-positive-solution} $0 < R(0) \leq K$, $\left( R(t) - R(0) \right) / t \to 0$, and therefore
			\begin{equation} \label{I-R-liminf}
				\liminf\limits_{t \to \infty} \frac{1}{t} \int\limits_{0}^{t} \left( I(s) + R(s) \right) ds
				= \frac{\theta + \kappa + \gamma}{\kappa + \gamma} \, \liminf\limits_{t \to \infty} \frac{1}{t} \int\limits_{0}^{t} I(s) \, ds,
			\end{equation}
			coming from the relation $(\kappa + \gamma) \int_{0}^{t} R(s) \, ds = \theta \int_{0}^{t} I(s) \, ds - R(t) + R(0)$, which follows immediately from the equation for $R$ in the system \eqref{sde1}.
			
			It remains to divide the expression \eqref{persistence-master-estimate} by $t$ and to let $t \to \infty$. First, note that according to the upper bound \eqref{logf-bounds} and $(\sigma_{2}^{\ast})$, it follows that
			\begin{equation*}
				\langle M^{G}, M^{G} \rangle_{t} = \int\limits_{0}^{t} \int\limits_{\mathbb{R}} \left( \frac{f^{\prime}\left( \beta(s-) \right)}{1 + f\left( \beta(s-) \right)} \right)^{2} \sigma_{s}^{2}\left( z, \beta(s-) \right) \Lambda(ds, dz)
				\leq \frac{K_{f^{\prime}}^{2}}{\left( 1 + f_{min} \right)^{2}} \left[ \left( g^{\ast}(0) \right)^{2} \Lambda^{B}_{t} + G_{2} \, \widehat{\Lambda}^{H}_{t} \right],
			\end{equation*}
			so that $M^{G}(t)/t \to 0$ $P$-a.s. by the SLLN recalled in \eqref{SLLN-consequence}.
			
			Further, $G \left( I(0), \beta(0) \right)/t \to 0$, the continuous part of the noise contributes by $\ell_{B}$ and the jumps contribute by $\ell_{H}$ by \eqref{jump-sum-limit}, both defined in \eqref{ell-definitions}. 
			
			Now taking the $\limsup_{t \to \infty}$ in the upper bound \eqref{persistence-master-estimate} divided by $t$, using \eqref{G-bounded-above} on the left and the relations $\limsup_{t} (u_{t} + v_{t}) \geq \limsup_{t} u_{t} + \liminf_{t} v_{t}$ and $\limsup_{t} (u_{t} + v_{t}) \leq \limsup_{t} u_{t} + \limsup_{t} v_{t}$ on the right, together with \eqref{I-R-liminf}, we arrive at
			\begin{align}
				0 \ & \geq \ \delta \widetilde{\beta}_{max} f_{min} - \left( \kappa_{1} + \theta \right) - \frac{K_{x f^{\prime}}}{1 + f_{min}} - C_{G} \left( \ell_{B} + \ell_{H} \right) \nonumber \\
				& \qquad - \frac{\delta \widetilde{\beta}_{max} \kappa_{1} f_{min}}{\kappa K} \cdot \frac{\theta + \kappa + \gamma}{\kappa + \gamma} \, \liminf\limits_{t \to \infty} \frac{1}{t} \int\limits_{0}^{t} I(s) \, ds, \quad P\text{-a.s.} \label{persistence-final-inequality}
			\end{align}
			Since $C_{G}$ is exactly the coefficient of $\ell_{B} + \ell_{H}$ in notation \eqref{persistence-threshold} in the statement of the theorem, the first four terms on the right-hand side add up to $\delta \widetilde{\beta}_{max} f_{min} - \Theta$. Therefore, from inequality \eqref{persistence-final-inequality} it follows that
			\			\begin{align*}
			\liminf\limits_{t \to \infty} \frac{1}{t} \int\limits_{0}^{t} I(s) \, ds
			\ & \geq \ \frac{\kappa K \left( \kappa + \gamma \right)}{\delta \widetilde{\beta}_{max} \kappa_{1} f_{min} \left( \theta + \kappa + \gamma \right)} \left( \delta \widetilde{\beta}_{max} f_{min} - \Theta \right) \\
			& = \ \frac{\kappa + \gamma}{\kappa + \gamma + \theta} \cdot \frac{\kappa K}{\kappa_{1}} \left( 1 - \frac{\Theta}{\delta \widetilde{\beta}_{max} f_{min}} \right) > 0, \quad P\text{-a.s.},
			\end{align*}
			which proves the theorem.
		\end{proof}

\begin{figure}[H]
	\centering
	\includegraphics[width=\textwidth]{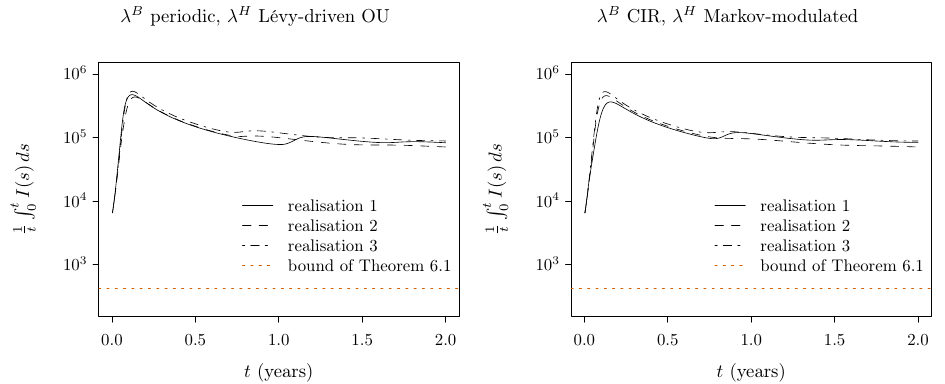}
	\caption{Persistence regime ($\underline{\mathcal{R}}_{0} = 1.16$ and $\overline{\mathcal{R}}_{0} = 3.32$ from Section~\ref{sec:simulation}): the time average $\frac{1}{t} \int_{0}^{t} I(s) \, ds$ for the three realisations of the noise from Figure~\ref{fig:beta} (solid, dashed and dash-dotted lines) and the two stochastic clocks, together with its lower bound from Theorem~\ref{thm-persistence}.} \label{fig:persistence-average}
\end{figure}

\begin{figure}[H]
	\centering
	\includegraphics[width=\textwidth]{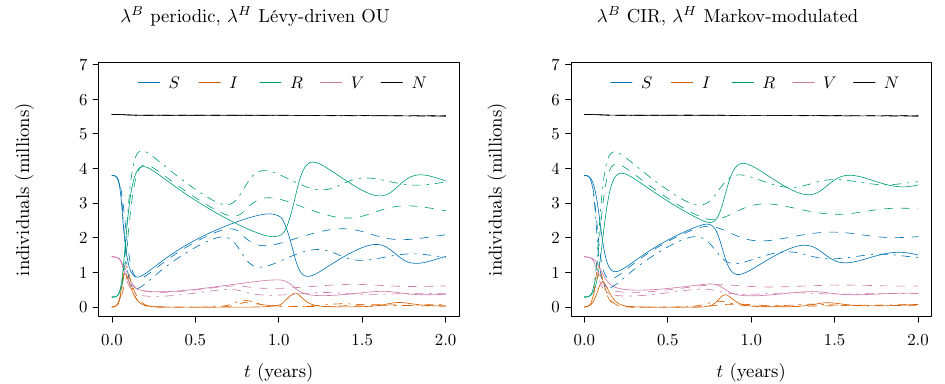}
	\caption{Persistence regime ($\underline{\mathcal{R}}_{0} = 1.16$ and $\overline{\mathcal{R}}_{0} = 3.32$ from Section~\ref{sec:simulation}): the compartments $S$, $I$, $R$, $V$ and the population size $N$ for the three realisations of the noise from Figure~\ref{fig:beta} (solid, dashed and dash-dotted lines) and the two stochastic clocks.} \label{fig:persistence-compartments}
\end{figure}

\begin{figure}[H]
	\centering
	\includegraphics[width=\textwidth]{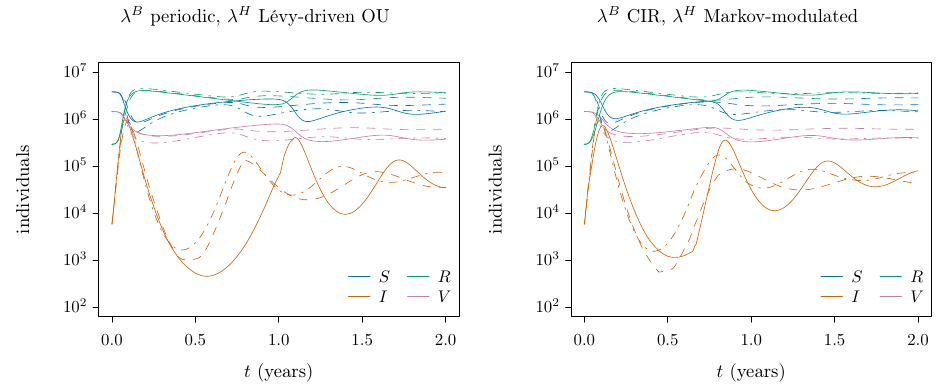}
	\caption{Persistence regime ($\underline{\mathcal{R}}_{0} = 1.16$ and $\overline{\mathcal{R}}_{0} = 3.32$ from Section~\ref{sec:simulation}): the compartments from Figure~\ref{fig:persistence-compartments} on a logarithmic scale.} \label{fig:persistence-log}
\end{figure}

\begin{figure}[H]
	\centering
	\includegraphics[width=\textwidth]{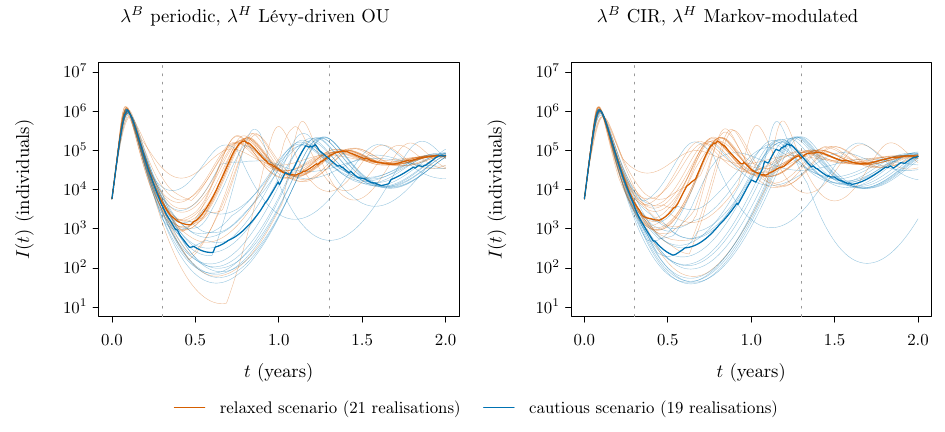}
	\caption{Persistence regime from Section~\ref{sec:simulation}: the number of infected individuals $I$ over two years for $40$ realisations of the noise and the two stochastic clocks, with the colours and the medians as in Figure~\ref{fig:ensemble-beta}.} \label{fig:ensemble-persistence}
\end{figure}

\begin{remark} \label{rem-persistence-interpretation}
	Assumption \eqref{persistence-condition} from Theorem \ref{thm-persistence}, under which the disease is persistent in mean, can be written in terms of the random basic reproduction number from Subsection \ref{subsec-stochastic-SIRV}, i.e. in terms of its lower bound $\underline{\mathcal{R}}_{0} = \delta \widetilde{\beta}_{max} f_{min} / (\kappa_{1} + \theta)$:
	\begin{equation*} \label{persistence-condition-R0}
		\underline{\mathcal{R}}_{0} \ > \ 1 + \frac{1}{\kappa_{1} + \theta} \left( \frac{K_{x f^{\prime}}}{1 + f_{min}} + \frac{K_{f^{\prime}}^{2} + K_{f^{\prime \prime}} \left( 1 + f_{min} \right)}{2 \left( 1 + f_{min} \right)^{2}} \left( \ell_{B} + \ell_{H} \right) \right) .
	\end{equation*}
	Since the terms added to $1$ are non-negative and $\mathcal{R}_{0}(t) \geq \underline{\mathcal{R}}_{0}$, condition \eqref{persistence-condition} implies that $\inf_{t \geq 0} \mathcal{R}_{0}(t) > 1$. In epidemiological terms, the disease persists if the random basic reproduction number $\mathcal{R}_{0}$ stays above $1$ even in the case in which the contact rate is at its minimal level and in which the vaccine doesn't provide a total protection against the disease, with a margin for the noise. The factor $\delta$ in $\underline{\mathcal{R}}_{0}$ comes from the bound $S + \delta V \geq \delta (S + V)$ -- for a perfectly protecting vaccine ($\delta = 0$) this bound reduces to $S \geq 0$ and gives no positive lower bound on the incidence. In that case $\underline{\mathcal{R}}_{0} = 0$ and condition \eqref{persistence-condition} cannot hold, which is why Theorem \ref{thm-persistence} requires the assumption $\delta > 0$. For the parameter values from Section~\ref{sec:simulation}, condition \eqref{persistence-condition} holds (see Table~\ref{tab:latent}), and Figure~\ref{fig:persistence-average} shows that the time average $\frac{1}{t} \int_{0}^{t} I(s) \, ds$ stays above its lower bound from Theorem~\ref{thm-persistence}, while Figures~\ref{fig:persistence-compartments}, \ref{fig:persistence-log} and~\ref{fig:ensemble-persistence} show the corresponding paths of the compartments and the paths of $I$ for $40$ realisations of the noise.
\end{remark}

\FloatBarrier
\section{Parameter values and simulation study}

\label{sec:simulation}
\suppressfloats[t]

In this section we present the parameter values and the simulated regimes of the random SIRV model \eqref{sde1}. The parameters of the model are of two different kinds. The epidemiological parameters in Table~\ref{tab:parameters} describe the disease and the population, and their values are set for the SARS-CoV-2 epidemic in Denmark, for which detailed public data and published estimates of all these quantities are available. Namely, the population size is taken from Statistics Denmark \cite{DST}, and the life expectancy, the vaccination rate and the reproduction number during the winter lockdown in 2020/2021 are derived from the Danish data collected by Our World in Data \cite{ArroyoMarioli2021,OWID, OWIDvacc}, while the generation time and the basic reproduction number before the first lockdown \cite{Valentin2021}, the infection fatality ratio \cite{Erikstrup2021, Pedersen2021}, the vaccine effectiveness and its waning \cite{Gram2022,Hansen2021medRxiv} and the protection after infection \cite{Hansen2021Lancet, Michlmayr2022} are taken from published Danish studies.

\begin{table}[ht]
	\centering
	\small
	\caption{Epidemiological parameters of the simulations (rates per year -- one day $= 1/365$ year).} \label{tab:parameters}
	\begin{tabular}{@{} >{\raggedright\arraybackslash}p{1.55cm} >{\raggedright\arraybackslash}p{3.3cm} >{\raggedright\arraybackslash}p{2.45cm} >{\raggedright\arraybackslash}p{7.4cm} @{}}
		\hline
		Symbol & Meaning & Value & Source \\
		\hline
		$K$ & carrying capacity & $5.85 \cdot 10^{6}$ & population of Denmark on 1 July 2021 \cite{DST} \\
		$\kappa$ & natural death rate & $0.0124$ & inverse life expectancy, $1/80.9$ years \cite{OWID} \\
		$\theta$ & recovery rate & $58.49$ & inverse mean generation time, $1/6.24$ days \cite{Valentin2021} \\
		$\kappa_{1}$ & death rate of the infected & $0.3063$ & $\kappa + \mathrm{IFR} \, \theta / (1 - \mathrm{IFR})$, with the infection fatality ratio $\mathrm{IFR} = 0.5\,\%$, lying between $0.089\,\%$ (ages 17--69) \cite{Erikstrup2021} and $5.4\,\%$ (ages $70+$) \cite{Pedersen2021} \\
		$\rho$ & vaccination rate & $1.095$ & $0.003$ per day, the median over 2021 of the daily first doses per unvaccinated person \cite{OWID, OWIDvacc} \\
		$\delta$ & susceptibility of $V$ relative to $S$ & $0.5$ & $1 - \mathrm{VE}$, where the vaccine effectiveness $\mathrm{VE}$ is the relative reduction of the risk of infection after vaccination, about $50\,\%$ \cite{Hansen2021medRxiv} \\
		$\alpha$ & waning of vaccine protection & $3.04$ & mean protection of $120$ days, in which the $\mathrm{VE}$ falls from $40\,\%$ to $12.6\,\%$ \cite{Gram2022} \\
		$\gamma$ & waning of immunity after infection & $1$ & mean protection of one year, between the Danish estimates \cite{Hansen2021Lancet, Michlmayr2022} \\
		$\widetilde{\beta}_{max}$ & ceiling of the contact rate & $195.22$ (persistence), $41.16$ (extinction) & $\overline{\mathcal{R}}_{0} (\kappa_{1} + \theta)$ by Remark~\ref{rem-R0}, with $\overline{\mathcal{R}}_{0} = 3.32$ before the first lockdown \cite{Valentin2021} and $\overline{\mathcal{R}}_{0} = 0.7$ during the winter lockdown in 2020/2021 \cite{OWID, ArroyoMarioli2021} \\
		\hline
	\end{tabular}
\end{table}

On the other hand, the underlying process $\beta$, the jump coefficient $\sigma$ and the stochastic clock $\Lambda$, i.e. the process $\lambda = (\lambda^{B}, \lambda^{H})$, are latent. Namely, the process $\beta$ enters the model only through the contact rate $\widetilde{\beta} = \widetilde{\beta}_{max} f(\beta)$, so that its influence is determined by the choice of the function $f$. The coefficient $\sigma$ describes the response of $\beta$ to shocks which are not observed, and the stochastic clock affects only the law of the noise $\mu$. Their estimation from incidence data is one of the open problems mentioned in Section~\ref{sec:conclusion}. Therefore, the values in Table~\ref{tab:latent} are not estimates, but rather chosen values so that the effects of the noise and of the stochastic clock are visible on the time scale of the epidemic, under the assumptions of Theorems~\ref{thm-extinction} and~\ref{thm-persistence}.

\begin{table}[!ht]
	\centering
	\small
	\caption{Parameters of the latent part of the model, chosen for the purpose of illustration; rates per year (one day $= 1/365$ year). The processes $\lambda^{B}$ and $\lambda^{H}$ of both pairs, from Example~\ref{ex-clocks}, all have $\ell^{B} = \ell^{H} = 1$ in $(l4)$.} \label{tab:latent}
	\begin{tabular}{@{} >{\raggedright\arraybackslash}p{1.9cm} >{\raggedright\arraybackslash}p{3.3cm} >{\raggedright\arraybackslash}p{2.6cm} >{\raggedright\arraybackslash}p{6.9cm} @{}}
		\hline
		Symbol & Meaning & Value & Choice \\
		\hline
		$f_{min}$ & lower bound of $f$ in \eqref{f-example} & $0.7$ & behaviour alone lowers the contact rate by at most $30\,\%$; Theorem~\ref{thm-persistence} needs $f_{min} > 1 / (\delta \overline{\mathcal{R}}_{0}) = 0.60$ \\
		$\beta(0)$ & initial value of $\beta$ & $0$ & neutral state, $\widetilde{\beta}(0) = 0.85 \, \widetilde{\beta}_{max}$ \\
		$s_{0}$, $s_{1}$ & coefficient \eqref{sigma-example-affine}, $h \equiv 1$ & $1.2$, $3$ & the long-run averages \eqref{ell-definitions} are $\ell_{B} = s_{0}^{2} \ell^{B} = 1.44$ and $\ell_{H} = s_{1}^{2} \ell^{H} \int_{\mathbb{R}_{0}} z^{2} \, \nu(dz) = 54.75$, and $\ell_{B} + \ell_{H} = 56.2$ is below the right-hand sides $254.7$ and $250.6$ of \eqref{extinction-condition} and \eqref{persistence-condition} \\
		$\nu$ & L\'{e}vy measure & $\nu(\mathbb{R}_{0}) = 6.08$ & compound Poisson, one shock per $60$ days of operational time; sizes $\mathrm{NIG}(1/\sqrt{2}, 0, 1/\sqrt{2}, 0)$ \\
		$\lambda^{B}$, pair 1 & periodic, Example~\ref{ex-clocks}(ii) & $c = 1$, $a = 0.5$, $p = 1$ & seasonal, with the period of one year \\
		$\lambda^{H}$, pair 1 & L\'{e}vy-driven Ornstein--Uhlenbeck, Example~\ref{ex-clocks}(v) & $a = 26.07$ & $1/a = 14$ days; the subordinator $L$ is compound Poisson with rate $2$ and exponential jumps of mean $0.5$ \\
		$\lambda^{B}$, pair 2 & Cox--Ingersoll--Ross, Example~\ref{ex-clocks}(iv) & $a = 26.07$, $b = 1$, $\varsigma = 5.73$ & $1/a = 14$ days \\
		$\lambda^{H}$, pair 2 & Markov-modulated, Example~\ref{ex-clocks}(iii) & $\lambda_{1} = 0.5$, $\lambda_{2} = 2$ & $J$ stays on average $60$ days in state $1$ and $30$ days in state $2$, $\pi = (2/3, 1/3)$ \\
		\hline
	\end{tabular}
\end{table}

The function $h$ in the state-dependent jump coefficient from Example~\ref{ex-f-sigma}(b), which describes how the response to a shock depends on the current value of $\beta$, is a modeling choice as well, and here we take $h = 1$, i.e. the coefficient \eqref{sigma-example-affine}. For this choice, from \eqref{beta-explicit} it follows that $|\beta(t)| \to \infty$, $P$-a.s., and since $f$ given by \eqref{f-example} is almost constant outside the interval $[-3, 3]$, the contact rate eventually settles at $\widetilde{\beta}_{max}$ or at $\widetilde{\beta}_{max} f_{min}$. The same behaviour was observed in simulations with several other choices for the bounded function $h$.

The stochastic clock is also a non-observable quantity, and the two pairs of processes $\lambda^{B}$ and $\lambda^{H}$ in Table~\ref{tab:latent}, taken from Example~\ref{ex-clocks}, are chosen for the purpose of illustration. In the first pair, the periodic intensity $\lambda^{B}$ describes a seasonal variation of the everyday fluctuations in contact behaviour, while the L\'{e}vy-driven Ornstein--Uhlenbeck intensity $\lambda^{H}$ describes short bursts of shock activity. In the second pair, the Cox--Ingersoll--Ross intensity $\lambda^{B}$ describes a randomly varying intensity of the fluctuations, while the Markov-modulated intensity $\lambda^{H}$ alternates between calm and turbulent periods, the latter with four times as many shocks (e.g. periods of frequent changes of the measures). All above mentioned examples of $\lambda^{B}$ and $\lambda^{H}$ have the same long-run average speed $\ell^{B} = \ell^{H} = 1$ in $(l4)$. Hence, the two stochastic clocks distribute the noise differently in time, but they lead to the same thresholds and bounds of Theorems~\ref{thm-extinction} and~\ref{thm-persistence}. In order to isolate the effect of the stochastic clock, the same paths of the processes $W$ and $N$ from Theorem~\ref{theorem-time-change} are used.

The simulation of the noise is based on Theorem~\ref{theorem-time-change}. Namely, the processes $W$ and $N$ are evaluated at the operational times of the stochastic clock, and the system \eqref{sde1} is then solved numerically in the R package \texttt{yuima} \cite{Iacus}. In the extinction regime ($\overline{\mathcal{R}}_{0} = 0.7$) the initial value is $(S(0), I(0), R(0), V(0)) = (0.55, 0.02, 0.15, 0.23) \, K$ and the epidemic is observed for half a year, while in the persistence regime ($\underline{\mathcal{R}}_{0} = 1.16$, $\overline{\mathcal{R}}_{0} = 3.32$) the initial value is taken to be $(0.65, 0.001, 0.05, 0.249) \, K$ and the epidemic is observed for two years. In both cases $N(0) = 0.95 \, K$, so that the initial value lies in the positively invariant set $\Gamma$ required by Theorems~\ref{thm-extinction} and~\ref{thm-persistence}, and the larger initial number of infected individuals in the extinction regime only makes the decay of $I$ more visible. In both regimes, $40$ realisations of the noise are simulated for each of the two stochastic clocks. The R script which derives the values in Table~\ref{tab:parameters} from the data and reproduces all simulations can be provided on request.

\FloatBarrier

\section{Concluding remarks} \label{sec:conclusion}

We have studied a random SIRV model in which the compartments satisfy ordinary differential equations and the only random object is the contact rate $\widetilde{\beta}$ given by the SDE \eqref{SDE-beta}. Contact rate $\widetilde{\beta}$ is a bounded transformation of the underlying process $\beta$ which is driven by the time-changed L\'{e}vy noise, i.e. by the martinagel random measure $\mu$. Therefore, the SIRV system is coupled with the SDE for the contact rate process $\widetilde{\beta}$, thus enabling the compartmental SIRV system to be treated pathwise. This approach resulted in a basic reproduction number of the random system which is itself a stochastic process, in particular just a linear transformation of the stochastic contact rate $\widetilde{\beta}$. 
	
The extinction and persistence thresholds derived in Theorems \ref{thm-extinction} and \ref{thm-persistence} can be expressed in terms of the random basic reproduction number $\mathcal{R}_{0}$ (i.e. in terms of its upper bound $\overline{\mathcal{R}}_{0}$ and its lower bound $\underline{\mathcal{R}}_{0}$, respectively), and the long-run averages $\ell_{B}$, $\ell_{H}$ of the noise which are of stochastic nature. The stochastic clock $\Lambda$ enters both thresholds only through $\ell_{B}$ and $\ell_{H}$, that is, through the long-run speed at which the operational time of the fluctuations and of the shocks accumulates. The simulations, with the parameter values from Section~\ref{sec:simulation} based on the data for SARS-CoV-2 in Denmark, illustrate both regimes and show that two stochastic clocks with the same long-run averages can lead to considerably different paths of the epidemic. 
	
A clock which runs faster on average makes both extinction and persistence thresholds harder to verify, a clock with a finite long-run average speed keeps the thresholds finite, and a clock which eventually stop would remove the stochastic terms $\ell_{B}$ and $\ell_{H}$ from the thresholds. Further, self-exciting clocks of Hawkes type, which would let the epidemic accelerate its own clock, are excluded by the conditional independence structure of the noise (Remark \ref{rem-Lambda-generalisation}). For models in which such self-excitation acts directly on the compartments we refer to \cite{Tymoshenko-2026}. 
	
Finally, this epidemic model raises two statistical problems, which we leave for future research. The first is the estimation of the parameters of the model, of the function $f$, and of the unobserved underlying process $\beta$, from observed incidence data. The second is the recovery of the unobserved stochastic clock $\Lambda$ from an observed path of the contact rate, a variant of the recovery problem for time-changed L\'{e}vy processes studied in \cite{MW}.

\bigskip \bigskip

\noindent{\textbf{Acknowledgements.}} 
	
\noindent This paper is based upon work from COST Action STOCHASTICA (CA24104), supported by COST (European Cooperation in Science and Technology), \url{https://www.ucc.ie/en/stochastica/}.
Giulia Di Nunno acknowledges support of the SURE-AI centre, RCN project nr. 357482.
Jasmina \DJ or\dj evi\'{c} is supported by the Ministry of Science, Technological Development and Innovation of Republic of Serbia, agreement no. 451-03-34/2026-03/200124. Nenad \v{S}uvak was supported by the Croatian Science Foundation project Scaling in Stochastic Models (HRZZ-IP-2022-10-8081). 
	
	\bigskip

	\begin{appendices}

	\section{Stochastic integration and quadratic variations} \label{app:noise}

	In this appendix we recall the construction of the stochastic integral with respect to the random measure $\mu$ given in Definition \ref{mu-definition}, which gives a meaning to the contact-rate model \eqref{SDE-beta}, and we collect the facts about quadratic variations used in the proofs. All the results are due to \cite{GS2}, while for the last paragraph, concerning general semimartingales, we refer to \cite{Protter}.

	We begin with the variations of the random measure $\mu$. Since it is a random measure and not a classical process indexed by time, its quadratic variation $[\mu, \mu]$ and predictable quadratic variation $\langle \mu, \mu \rangle$ are again measures on $\mathcal{B}_{\mathcal{X}}$, given by
	\begin{align} \label{mu-variations}
		[\mu, \mu](dt, dz) & = \mathbf{1}_{\{ z = 0 \}} \lambda^{B}_{t} \, dt + \mathbf{1}_{\{ z \neq 0 \}} H(dt, dz), \nonumber \\
		\langle \mu, \mu \rangle(dt, dz) & = \mathbf{1}_{\{ z = 0 \}} \lambda^{B}_{t} \, dt + \mathbf{1}_{\{ z \neq 0 \}} \lambda^{H}_{t} \nu(dz) \, dt = \Lambda(dt, dz),
	\end{align}
	and they differ only in the jump part, where the predictable intensity $\lambda^{H}_{t} \nu(dz) \, dt$ of the shocks in $\langle \mu, \mu \rangle$ is replaced by their actual count $H(dt, dz)$ in $[\mu, \mu]$. Note that the predictable quadratic variation of the random measure $\mu$ is exactly the random measure $\Lambda$, determining the operational time in which $\mu$ runs. This is the precise sense in which $\Lambda$ is called the stochastic clock of $\mu$.

	Both quadratic and predictable quadratic variations in \eqref{mu-variations} are finite on every set $\Delta$ with $\Lambda(\Delta) < \infty$ $P$-a.s., which by (l3) covers $\Delta = [0, t] \times \{ 0 \}$ and $\Delta = [0, t] \times A$ with $\nu(A) < \infty$, for every $t > 0$. However, they are finite on the whole $[0, t] \times \mathbb{R}$ if and only if the L\'{e}vy measure $\nu$ is finite (i.e., $\nu(\mathbb{R}_{0}) < \infty$), in which case
	\begin{equation} \label{mu-variations-finite-activity}
		[\mu, \mu] \left( [0, t] \times \mathbb{R} \right) = \Lambda^{B}_{t} + H \left( [0, t] \times \mathbb{R}_{0} \right), \qquad
		\langle \mu, \mu \rangle \left( [0, t] \times \mathbb{R} \right) = \Lambda^{B}_{t} + \nu(\mathbb{R}_{0}) \widehat{\Lambda}^{H}_{t},
	\end{equation}
	both finite $P$-a.s. by $\it{(l3)}$, the first one because a Poisson random measure with finite L\'{e}vy measure produces only finitely many jumps on a bounded time-interval. Since the finite activity is not assumed in this paper (see Remark \ref{rem-finite-activity}), the expressions \eqref{mu-variations-finite-activity} are not available in general, and the quantities that remain finite are the variations of the integral process introduced below.

	According to \cite{GS2}, Section $2$, the stochastic integral with respect to the random measure $\mu$ is defined for $\mathbb{F}$-predictable random fields $\varphi$ satisfying
	\begin{equation} \label{integrand-condition}
		E \left[ \int\limits_{0}^{t} \int\limits_{\mathbb{R}} \varphi_{s}^{2}(z) \, \Lambda(ds, dz) \right] < \infty, \qquad \text{for every } t > 0,
	\end{equation}
	as the $L^{2}$-limit of the integrals of simple random fields, and it satisfies the isometry
	\begin{equation} \label{mu-isometry}
		E \left[ \left( \int\limits_{0}^{t} \int\limits_{\mathbb{R}} \varphi_{s}(z) \, \mu(ds, dz) \right)^{2} \right] = E \left[ \int\limits_{0}^{t} \int\limits_{\mathbb{R}} \varphi_{s}^{2}(z) \, \Lambda(ds, dz) \right], \qquad \text{for every } t > 0.
	\end{equation}
The integral process
	\begin{equation*} \label{M-integral-process}
		M(t) = \int\limits_{0}^{t} \int\limits_{\mathbb{R}} \varphi_{s}(z) \, \mu(ds, dz), \qquad t \geq 0,
	\end{equation*}
	is a square-integrable $\mathbb{F}$-martingale: the second moment is finite by the isometry \eqref{mu-isometry}, while the martingale property is inherited from the driving noise $\mu$, established in Section \ref{subsec-preliminaries}, together with the $\mathbb{F}$-predictability of the random field $\varphi$. By localisation, the construction extends to random fields $\varphi$ with $\int_{0}^{t} \int_{\mathbb{R}} \varphi_{s}^{2}(z) \, \Lambda(ds, dz) < \infty$ $P$-a.s. for every $t > 0$, in which case $M$ is a locally square-integrable martingale. Its quadratic variation and its predictable quadratic variation are obtained by integrating $\varphi^{2}$ against the corresponding variations \eqref{mu-variations} of the random measure $\mu$:
	\begin{equation*} \label{mu-integral-qv}
		[M, M]_{t} = \int\limits_{0}^{t} \varphi_{s}^{2}(0) \lambda_{s}^{B} \, ds + \int\limits_{0}^{t} \int\limits_{\mathbb{R}_{0}} \varphi_{s}^{2}(z) \, H(ds, dz),
		\qquad
		\langle M, M \rangle_{t} = \int\limits_{0}^{t} \int\limits_{\mathbb{R}} \varphi_{s}^{2}(z) \, \Lambda(ds, dz),
	\end{equation*}
	both finite $P$-a.s. for every $t > 0$ (they have equal expectations, and $E \left[ \langle M, M \rangle_{t} \right]< \infty$ by \eqref{integrand-condition}). Note that, unlike \eqref{mu-variations-finite-activity}, the finiteness of $[M, M]_{t}$ and $\langle M, M \rangle_{t}$ requires no finiteness of the L\'{e}vy measure $\nu$. Further, note that the jumps of $M$ are $\Delta M(s) = \int_{\mathbb{R}_{0}} \varphi_{s}(z) \, H(\{ s \}, dz)$, so that the second integral in $[M, M]_{t}$ is the sum $\sum_{0 < s \leq t} \left( \Delta M(s) \right)^{2}$ of the squared jumps, which converges even when $M$ has infinitely many small jumps on bounded time-intervals.

	Finally, we recall two facts about quadratic variations of general semimartingales, used in the proofs in Sections \ref{sec:EUPsolution}, \ref{sec:extinction} and \ref{sec:persistence} (see \cite{Protter}, Chapters II and III). For a semimartingale $X$, the quadratic variation $[X, X]$ is a c\`{a}dl\`{a}g, non-decreasing, adapted process such that $\Delta [X, X]_{t} = (\Delta X_{t})^{2}$, and it decomposes path by path into its continuous part and its jump part,
	\begin{equation} \label{relation_quad_pred}
		[X, X]_{t} = [X, X]^{c}_{t} + \sum\limits_{0 < s \leq t} (\Delta X_{s})^{2} ,
	\end{equation}
	where $[X, X]^{c} = \langle X^{c}, X^{c} \rangle$ is the quadratic variation of the continuous local martingale part $X^{c}$ of $X$. If $[X, X]$ is locally integrable, its compensator $\langle X, X \rangle$, called the predictable quadratic variation of $X$, is the unique predictable non-decreasing process such that $[X, X] - \langle X, X \rangle$ is a local martingale. For the process $\beta$, by whose transformation we define the contact rate $\widetilde{\beta}$, both quadratic and predictable quadratic variation are computed in \eqref{beta-qv}, and it is the decomposition \eqref{relation_quad_pred} that in the It\^{o} formula separates the contribution of the continuous fluctuations from the contribution of the jumps.
\end{appendices}

\end{document}